\documentclass[mathpazo]{cicp}

\usepackage[utf8]{inputenc}
\usepackage{amssymb,amsmath,bm}
\usepackage{amsthm}
\usepackage{mathrsfs}
\usepackage{float}
\usepackage{cite}
\usepackage[colorlinks,linkcolor=blue,anchorcolor=blue,citecolor=red]{hyperref}
\usepackage{cleveref}
\usepackage{graphicx}
\usepackage{hyperref}
\usepackage{subfigure}
\numberwithin{equation}{section}
\numberwithin{figure}{section}
\usepackage{subcaption}

\usepackage{color} 
\usepackage{diagbox}

\usepackage{verbatim}

\usepackage{url}
\usepackage[colorlinks,linkcolor=blue,anchorcolor=blue,citecolor=red]{hyperref}

\usepackage{algpseudocode,algorithm,algorithmicx}
\usepackage{algcompatible}

\usepackage{tabu}                     
\usepackage{multirow}                 
\usepackage{multicol}                 
\usepackage{multirow}                
\usepackage{float}                    
\usepackage{makecell}                 
\usepackage{booktabs}                 
\def\x{\boldsymbol{x}}

\def\rh{\boldsymbol{\rho}}

\def\v{\mathbf{v}}
\def\V{\mathbf{V}}

\def\dd{\mathbf{d}}

\def\f{\mathbf{f}}

\def\h{\mathbf{h}}

\def\m{\mathbf{m}}

\def\n{\mathbf{n}}
\def\u{\mathbf{u}}

\def\t{\mathrm{t}}

\def\d{\mathrm{d}}
\def\eff{\mathrm{eff}}
\def\div{\mathrm{div}}

\def\eps{\epsilon}

\def\0{\mathbf{0}}

\def\<{\langle}
\def\>{\rangle}

\definecolor{lzcol}{rgb}{1, 0, 0}

\definecolor{drcol}{rgb}{0,0,1}

\definecolor{jrcol}{rgb}{0,1,0}

\DeclareMathOperator*{\argmin}{arg\,min}

\begin{document}
\title{Numerical Homogenization of Landau-Lifshitz Equation with Rough Coefficients}


\author[Zetao Ma et.~al.]{Zetao Ma\affil{1}, Jingrun Chen\affil{2,3}, Rui Du\affil{4,5}\comma\corrauth ,  and  Lei Zhang\affil{1,6,7}}

\address{\affilnum{1}\ School of Mathematical Sciences,
         Shanghai Jiao Tong University, Shanghai 200240, P.R. China. \\
           \affilnum{2}\ School of Mathematical Sciences, University of Science and Technology of China, Hefei, Anhui 230026, P.R. China.\\
           \affilnum{3}\ Suzhou Institute for Advanced Research, University of Science and Technology of China and Suzhou Big Data \& AI Research and Engineering Center, Suzhou, Jiangsu 215123, P.R. China. \\
           \affilnum{4}\ School of Mathematical Sciences, Soochow University, Suzhou, 215006, P.R. China. \\
           \affilnum{5}\ Mathematical Center for Interdisciplinary Research, Soochow University, Suzhou, 215006, P.R. China. \\
           \affilnum{6}\ Institute of Natural Sciences, Shanghai Jiao Tong University, Shanghai 200240, P.R. China. \\
           \affilnum{7}\ MOE-LSC, Shanghai Jiao Tong University, Shanghai 200240, P.R. China. \\
           }
    \emails{{\tt 770120068@qq.com} (Z.~Ma), {\tt jingrunchen@ustc.edu.cn} (J.~Chen), {\tt durui@suda.edu.cn} (R.~Du), and {\tt  lzhang2012@sjtu.edu.cn} (L.~Zhang)}



\begin{abstract}
In this work, we develop a numerical homogenization approach for the fully nonlinear Landau-Lifshitz equation with rough coefficients, including non-periodicity and nonseparable scales. Direct numerical resolution of such multiscale problems on fine meshes incurs prohibitive computational costs. To address this challenge, we propose an efficient coarse scale approximation through localized basis functions derived from energy minimization within the Generalized Rough Polyharmonic Splines (GRPS) framework. These basis functions preserve critical multiscale features while operating on a computationally tractable coarse mesh. The nonlinear, vectorial, and non-symmetric nature of the Landau-Lifshitz equation necessitates careful design of variational formulations for basis construction. We introduce several such formulations, each tailored to specific structural aspects of the problem. Through systematic numerical experiments, we demonstrate that our approach achieves significant computational savings without compromising accuracy, offering a robust framework for simulating multiscale magnetic systems with complex microstructures.
\end{abstract}


\maketitle

\section{Introduction}
\label{sec1}

The magnetization dynamics in a ferromagnetic material is characterized by the Landau-Lifshitz (LL) equation \cite{landau1992theory, gilbert1955lagrangian} for the magnetization \(\m: \Omega \times (0, T] \rightarrow \mathbb{R}^3\), which assumes the following dimensionless form,
\begin{equation}\label{eqn:LL}
\left\{
    \begin{array}{ccc}  
        \partial_t \m = -\m \times \h_{\eff} - \lambda \m \times (\m \times \h_{\eff}),  & \x \in \Omega, \, \t \in (0,T], \\
        \kappa(\x) \frac{\partial \m}{\partial \n} = \0,  & \x \in \partial \Omega, \, \t \in (0,T], \\
        \m(\x,0) = \m_0(\x),  & \x \in \Omega,
    \end{array}
    \right.
\end{equation}
where \(\Omega \subseteq \mathbb{R}^d\), \(d=2,3\), is a bounded, smooth, and convex spatial domain. The parameter \(\lambda > 0\) represents the dimensionless damping parameter, indicating the magnitude of the damping effect, and the initial data \(\m_0\) satisfies \(|\m_0| = 1\) in the point-wise sense. The vector \(\n\) represents the outward unit normal to the boundary \(\partial \Omega\). The effective field \(\h_{\eff}\) combines the exchange and anisotropic terms, and is given by
\begin{equation}\label{effectivefield}
    \h_{\eff}[\m] = \underbrace{\div(\kappa \nabla \m)}_{\text{exchange\; term}} - \underbrace{(\m - (\m \cdot \u)\u)}_{\text{anisotropic\; term}},
\end{equation}
where the exchange parameter \(\kappa(\x)\) consists of rough coefficients (i.e., \(\kappa \in L^{\infty}\)), such that
\begin{equation}\label{upper_lower_bound}
    \kappa_{\min} |\xi|^2 \le \xi^T \kappa(\x) \xi \le \kappa_{\max} |\xi|^2,
\end{equation}
for all \(\xi \in \mathbb{R}^d\) and almost every \(\xi \in \Omega\). The material dependence of the exchange parameters varies at a very small scale and may not rely on concepts of ergodicity or scale separation, such as in high performance rare-earth magnetic alloys \cite{gutfleisch2011magnetic}. Studying these materials is challenging due to their strong heterogeneity and highly nonlinear nature. Additionally, we denote the anisotropy energy density as 
\begin{equation}
\m_a = \m - (\m \cdot \u)\u = (0, m_2, m_3)^T
\label{eqn:anisotropyenergydensity}
\end{equation}
with \(\m = (m_1, m_2, m_3)^T\) and  \(\u = (1,0,0)^T\) (i.e., along the X-axis) for a uniaxial material. The Landau-Lifshitz energy functional is defined as
\begin{equation}\label{LLenergy}
    F[\m] = \frac{1}{2} \int_{\Omega} \left( \kappa |\nabla \m|^2 + (|\m|^2 - |\m \cdot \u|^2) \right) \d \x = -\frac{1}{2}(\h_{\eff}[\m], \m),
\end{equation}
where the first equality indicates that $F[\m]$ is positive, and due to the minus sign ensures the non-positivity of $\h_{\eff}$.

Several papers have analyzed the homogenization of the LL model from the analytical perspective, where $\kappa(x)$ can be expressed as $\kappa(x) = a\left(\frac{x}{\varepsilon}\right)$, $a\left(x,\frac{x}{\varepsilon}\right)$, or $a\left(\frac{x}{\varepsilon}, \omega\right)$, showcasing various scale separation characteristics \cite{santugini2007homogenization, alouges2015homogenization, choquet2018homogenization, alouges2021stochastic, leitenmaier2022homogenization, chen2022multiscale, leitenmaier_Phdthesis}. In cases with scale separation ($\varepsilon \ll 1$) such that $\kappa(x) = a^\varepsilon$, the heterogeneous multiscale method (HMM) computes the homogenized solution using microscopic cell problems, reducing the computational burden associated with resolving fine scales \cite{UpscalingHMM, leitenmaier2022heterogeneous, leitenmaier2023finite}. However, these models overlook the exchange term with rough coefficients, which is a more practical scenario \cite{gutfleisch2011magnetic}, and cannot be effectively addressed by numerical methods based on the scale separation assumption.

To date, numerical homogenization without scale separation has matured for benchmark problems such as elliptic PDEs. Key methods include asymptotic homogenization \cite{papanicolau1978asymptotic,jikov2012homogenization}, numerical homogenization \cite{dur91,ab05,weh02}, heterogeneous multiscale methods \cite{ee03,abdulle2014analysis,ming2005analysis}, multiscale network approximations \cite{berlyand2013introduction}, multiscale finite element methods (MsFEM) \cite{Arbogast_two_scale_04,egw10,eh09}, variational multiscale methods \cite{hughes98,bazilevs2007variational}, flux norm homogenization \cite{berlyand2010flux,Owhadi:2011}, rough polyharmonic splines (RPS) and its generalization (GRPS) \cite{zhang_RPS,owhadi2017multigrid,liuGRPS}, localized orthogonal decomposition (LOD) \cite{MalPet:2014,altmann2021numerical,hauck2023super}, and generalized multiscale finite element methods (GMsFEM) \cite{egh12,chung2014adaptiveDG,chung2015residual,chung2014adaptive,ma2022novel}.

In this work, we employ RPS/GRPS as our primary numerical homogenization approach. These methods have proven effective for obtaining accurate multiscale solutions to elliptic, parabolic, and hyperbolic equations \cite{Owhadi:2011,parabolic_homo_zhang,zhang_RPS,owhadi2017gamblets,liuGRPS}. The GRPS method, as developed in \cite{liuGRPS}, incorporates multiple energy forms and measurement functions (edge, volume, and derivative measurements). Specifically, the edge-based GRPS (GRPS-E) achieves first-order accuracy in the $H^1$ sense for $g \in L^2(\Omega)$, while both volume-based (GRPS-V) and derivative-based (GRPS-D) approaches attain second-order accuracy when $g \in H^1(\Omega)$.

The Landau-Lifshitz equation \eqref{eqn:LL} presents several unique challenges beyond standard linear PDEs:

\begin{itemize}
\item \textbf{Nonlinearity:} The terms $-\m \times \h_{\eff} - \lambda \m \times (\m \times \h_{\eff})$ exhibit third-order nonlinearity. Various temporal discretization approaches \cite{cimrak2007survey,kruzik2006recent,GarcaCervera2007NUMERICALMA,AnewschemeLLG,gao2014optimal,praetorius2018convergence,di2020linear,xie2020second} have been developed, ranging from explicit (with strict stability constraints) to fully implicit (requiring expensive nonlinear solves) and semi-implicit schemes (offering better stability-efficiency balance). These choices directly affect the resulting bilinear forms in numerical homogenization.

\item \textbf{Skew-symmetry:} The term $-\m \times \h_{\eff}$ introduces a skew-symmetric component to the bilinear form, necessitating a generalized Lax-Milgram framework \cite{chen1998second}.

\item \textbf{Vectorial nature:} As a 3-component vectorial PDE, the LL equation triples the system dimensionality. The bilinear form can be constructed either as a combined variational form or as separate forms for each component.

\item \textbf{Constraint preservation:} The magnetization constraint $|\m| = 1$ is typically maintained through point-wise projection, though well-designed schemes can preserve this property without explicit projection.
\end{itemize}

While numerical homogenization has been successfully applied to nonlinear PDEs \cite{henning2020computational,henning2022superconvergence,maier2022multiscale,liu2021iterated,verfurth2022numerical}, vectorial problems (e.g., H(curl)-systems) \cite{gallistl2018numerical}, and non-symmetric problems (e.g., convection-dominated systems) \cite{freese2024computational,bonizzoni2024reduced,bonizzoni2022super}, this work represents the first application to the LL equation incorporating all these features simultaneously.

Using GRPS, we construct low-dimensional, operator-adapted basis functions with exponential decay, though LOD and GMsFEM remain viable alternatives. Through systematic investigation of different GRPS formulations (including bilinear forms, measurement functions, time discretizations, and anisotropy terms), we identify optimal configurations exhibiting varying convergence rates, including higher-order convergence.

The paper is organized as follows: Section \ref{sec:numericalmethod} overviews FEM approaches for the LL equation; Section \ref{sec3} details GRPS basis derivation and reduced-space approximation; Section \ref{sec4} presents numerical experiments validating our method's accuracy and efficiency; and Section \ref{sec:conclusion} provides concluding remarks.

\paragraph{Notation:} $(\cdot,\cdot)$ represents $L^2$-inner product in space.

\section{Numerical Methods for Model Problem}
\label{sec:numericalmethod}

The LL equation has been extensively studied using conventional numerical methods such as finite difference and finite element approaches. For comprehensive reviews of these standard techniques, we refer readers to \cite{cimrak2007survey,kruzik2006recent,GarcaCervera2007NUMERICALMA}. 

In our work, we focus on semi-implicit time discretization schemes \cite{AnewschemeLLG,gao2014optimal,praetorius2018convergence,di2020linear,xie2020second,cai2022second}, which offer a practical balance between numerical stability and computational efficiency. We consider a uniform partition of the time interval $[0, T]$ denoted by $0 = t_0 < t_1 < \cdots < t_N = T$, where the time step size is given by $\Delta t = T/N$.
At each time step, these schemes require solving for the updated magnetization field $\m^{n+1}\in \V := [H^1(\Omega)]^3$ through the variational formulation:

\begin{equation}\label{eqn:abstract_blinear_form}
    A^n(\m^{n+1}, \v) = (\f^n,\v), \quad \forall \v \in \V.
\end{equation}

The bilinear form $A^n$ can be decomposed into two components that will be specified later,
\begin{itemize}
    \item $B^n(\m^{n+1}, \v)$: Contains the essential multiscale structure and heterogeneity of the problem, which plays a crucial role in defining our coarse approximation space $V_H$ (where $H$ represents the coarse scale parameter);
    \item $C^n(\m^{n+1}, \v)$: Handles the linearized approximation of nonlinear terms, which as noted in \cite{doding2022two}, often has secondary importance when constructing the coarse space.
\end{itemize}

The presence of rough coefficients $\kappa(\x)$ with fine-scale (characterized by length scale $0<\eps\ll 1$) variations makes direct numerical simulation computationally prohibitive using conventional methods. This motivates our primary focus on developing an effective coarse space construction that can capture the essential multiscale behavior without resolving all fine-scale details.

\subsection{Time Discretization Schemes}
\label{sec:numericalmethod:timediscretization}

Using the fundamental constraint $|\m| = 1$ of micromagnetics, we can reformulate the LL equation \eqref{eqn:LL} into an equivalent form that is often more amenable to discover the bilinear form:

\begin{equation}\label{LLform1}
    \m_t - \lambda \h_{\eff} + \m \times \h_{\eff} = -\lambda (\m \cdot \h_{\eff})\m
\end{equation}

Several semi-implicit backward Euler-type schemes have been developed to handle the nonlinear terms while maintaining $O(\Delta t + h^2)$ accuracy in the $L^2$ norm. These schemes all share a common structure for the $B^n$ component of the bilinear form:

\begin{equation}\label{eqn:B_blinear_form}
    B^n(\m^{n+1},\v) := \lambda (\kappa \nabla \m^{n+1}, \nabla \v) + \lambda (\m_a^{n+1},\v) - (\m^n \times \kappa \nabla \m^{n+1}, \nabla \v)
\end{equation}

while differing in their treatment of the remaining terms:

\begin{itemize}
    \item \textbf{Cimr\'ak's scheme} \cite{cimrak2007survey}:
    \begin{align*}
        C^n_{Cimrak}(\m^{n+1},\v) &= \frac{1}{\Delta t}(\m^{n+1},\v) - (\m^n \times \m_a^{n+1},\v) + \lambda(\m^n \cdot \h_{\eff}^{n+1},\v), \\
        (\f^n_{Cimrak},\v) &= \frac{1}{\Delta t}(\m^n,\v);
    \end{align*}
    
    \item \textbf{Gao's scheme} \cite{gao2014optimal}:
    \begin{align*}
        C^n_{Gao}(\m^{n+1},\v) &= \frac{1}{\Delta t}(\m^{n+1},\v) - (\m^n \times \m_a^{n+1},\v), \\
        (\f^n_{Gao},\v) &= -\lambda((\m^n \cdot \h_{\eff}^n)\m^n,\v) + \frac{1}{\Delta t}(\m^{n},\v);
    \end{align*}
    
    \item \textbf{An's scheme} \cite{anrongbackward}:
    \begin{align*}
        C^n_{An}(\m^{n+1},\v) &= \frac{1}{\Delta t}(\m^{n+1},\v) - (\m^n \times \m_a^{n+1},\v) + \lambda ((\m^{n} \cdot \Bar{\h}_{\eff}^{n+1}) \cdot \m^{n}, \v), \\
        (\f^n_{An},\v) &= \frac{1}{\Delta t}(\m^{n},\v).
    \end{align*}
\end{itemize}

\subsection{Numerical Homogenization Framework}
\label{sec:prelim:numhom}

To construct our coarse approximation space $V_H$, we employ the Generalized Rough Polyharmonic Splines (GRPS) method \cite{zhang_RPS,liuGRPS,owhadi2017multigrid}, which has proven effective for numerical homogenization of multiscale problems. The approach begins with the following setup:

\begin{itemize}
    \item \textbf{Mesh structure}: Let $\mathcal{T}_H$ be a coarse simplicial subdivision of the domain $\Omega$ into $2N_c^2$ parts with maximum element diameter $H$. $\mathcal{T}_H$ can be uniformly refined $J$ times to create a fine mesh $\mathcal{T}_h$ with $h = 2^{-J}H$.
    
    \item \textbf{Function spaces}: On $\mathcal{T}_h$, we define the standard first-order conforming finite element space:
    \begin{equation*}
        V_h = \{ \varphi \in C^0(\Omega) \, | \, \varphi|_\tau \in \mathcal{P}_1(\tau), \, \forall \tau \in \mathcal{T}_h \} \subset H^1(\Omega)
    \end{equation*}
    with the vector-valued counterpart $\mathbf{V}_h = V_h^3$ for our LL equation application.
\end{itemize}

Consider the model elliptic problem $-\nabla \cdot (\kappa(\x) \nabla u) = f$ where $\kappa \in L^\infty(\Omega)$ is bounded and positive definite. The associated bilinear form is:

\begin{equation}
    B(u,v) = \int_\Omega \kappa(\x) \nabla u \cdot \nabla v \, \mathrm{d}x.
\end{equation}

The numerical homogenization seeks a coarse space $V_H \subset V_h$ of dimension much smaller than $V_h$ such that the coarse solution $u_H \in V_H$ satisfies:

\begin{equation}\label{eqn:numhom}
    B(u_H, v_H) = (f, v_H), \quad \forall v_H \in V_H
\end{equation}
with an approximation error comparable to the coarse mesh size $H$.

\subsubsection{GRPS Basis Construction}

The GRPS approach constructs basis functions through the following procedure \cite{liuGRPS}:

\begin{enumerate}
    \item Define $N_H$ measurement functions $\{\phi_i\}_{i=1}^{N_H}$ (typically characteristic functions of coarse patches or edges);
    
    \item For each $\phi_i$, define the constrained space:
    \begin{equation*}
        \mathcal{V}_i = \{\psi \in V_h \, | \, (\psi, \phi_j) = \delta_{ij}, \, j=1,...,N_H\};
    \end{equation*}
    
    \item Solve the constrained minimization problem for each basis function:
    \begin{equation}\label{eqn:phi}
        \psi_i = \argmin_{\psi \in \mathcal{V}_i} \|\psi\|_{B,\Omega}^2,
    \end{equation}
    where $\|\psi\|_{B,\Omega} = \sqrt{B(\psi,\psi)}$ is the energy norm.
\end{enumerate}

\begin{remark}
The choice of measurement functions significantly impacts both accuracy and computational efficiency. We focus on
\begin{itemize}
    \item \textbf{Volume-based (GRPS-V)}: Uses characteristic functions of coarse elements, with the set of measurement functions denoted as $\Phi_{\mathcal{T}}$;
    \item \textbf{Edge-based (GRPS-E)}: Uses characteristic functions of coarse edges, with the set of measurement functions denoted as $\Phi_{\mathcal{E}}$.
\end{itemize}
They provide optimal balance between accuracy ($O(H)$ for $f\in L^2(\Omega)$) and reasonable number of degrees of freedom \cite{liuGRPS}, while for GRPS-V we achieve higher accuracy $O(H^2)$ for $f \in H^1(\Omega)$.
\end{remark}

\subsubsection{Localization and Approximation Properties}

While the exact GRPS basis functions $\{\psi_i\}$ possess global support, their exponential decay property permits effective localization through the following construction. First, let $\Omega_i^0$ denote the nodal support of $\psi_i$. We then define the $\ell$-layer patch $\Omega_i^\ell$ as the $\ell$-th order neighborhood of $\Omega_i^0$. The localized basis function $\psi_i^\ell$ is obtained by solving the constrained minimization problem
\begin{equation}\label{eqn:philoc}
\psi_i^\ell = \argmin_{\substack{\psi \in \mathcal{V}_i \\ \mathrm{supp}(\psi)\subset\Omega_i^\ell}} \|\psi\|_{B,\Omega_i^\ell}^2,
\end{equation}
where the minimization is performed over the subspace $\mathcal{V}_i$ with support restricted to $\Omega_i^\ell$, and $B$ represents the appropriate energy norm for the problem.

\begin{theorem}[Localization Error \cite{zhang_RPS,liuGRPS}]\label{theorem_convergence_rate}
For sufficiently large $\ell \geq C\log(1/H)$, the localized basis satisfies:
\begin{enumerate}
    \item Exponential decay: 
    \begin{equation*}
        \|\psi_i - \psi_i^\ell\|_{B,\Omega} \leq C_1 e^{-C_2\ell} \|\psi_i\|_{B,\Omega};
    \end{equation*}
    
    \item Approximation quality:
    \begin{equation*}
        \|u - u_H^\ell\|_{B,\Omega} \leq C(H^{s+1} + e^{-C_2\ell})\|f\|_{H^s(\Omega)}
    \end{equation*}
    for $s=0,1$, where $u_H^\ell$ is the solution in $\mathrm{span}\{\psi_i^\ell\}$.
\end{enumerate}
\end{theorem}

\begin{remark}
Recent developments in super-localization techniques \cite{hauck2023super} and LOD/SLOD methods \cite{bonizzoni2024reduced,bonizzoni2022super} offer promising alternatives, particularly for problems with additional structural challenges like strong convection or skew-symmetric terms. While we focus on GRPS for this work, these approaches present interesting directions for future research.
\end{remark}

\section{GRPS for the LL equation}
\label{sec3}

In this section, we construct a coarse space for the bilinear form \eqref{eqn:abstract_blinear_form} by developing a multiscale basis for the magnetization components \(m_i\) (\(i=1,2,3\)). Solving \eqref{eqn:abstract_blinear_form} within this framework allows us to accurately capture the multiscale nature of the problem. This task is particularly challenging for the LL equation, as the bilinear form \(B^n(\m^{n+1}, \v)\) in \eqref{eqn:B_blinear_form} does not conform to a simple elliptic structure. Additionally, the choice of measurement function \cite{liuGRPS} plays a crucial role, as it can significantly influence convergence rates.

\subsection{Construction of the Multiscale Basis Functions}

The multiscale basis functions are constructed by solving the constrained minimization problem \eqref{eqn:phi}, representing the elliptic problem's energy norm. For the LL equation, the bilinear form \( B^n(\m, \v) \) in \eqref{eqn:B_blinear_form} induces an energy norm equivalent (modulo $\lambda$) to the LL energy functional:
\[
F[\m] = \frac{1}{2} \int_{\Omega} \left( \kappa |\nabla \m|^2 + (|\m|^2 - |\m \cdot \u|^2) \right) \d\x = -\frac{1}{2}(\h_{\eff}[\m], \m),
\]
as defined in \eqref{LLenergy}.

Given measurement functions \(\{\phi_j\}_{j=1}^{N_H} \subset \Phi\) ($\Phi = \Phi_{\mathcal{E}}$ or $\Phi_{\mathcal{T}}$), the GRPS basis solves:
\[
\psi_i = \underset{\m}{\text{argmin}} \, F[\m] \quad \text{s.t.} \quad \langle \phi_j, \m \rangle = \delta_{ij}, \quad \forall j.
\]
The resulting GRPS space $V_H = \text{span}\{\psi_i\}_{i=1}^{N_H}$ replaces the fine-scale space $V_h$ on coarse meshes, and can apply to time discretization schemes such as  \eqref{Cimark}, \eqref{GaoFEM}, \eqref{AnFEM}. 

To obtain the local GRPS basis space $\Psi^l = \text{span}\{\psi_i^l\}$ within the local subdomain $\Omega_i^l$, we solve the following optimization problem:
\begin{equation}\label{vector_optimize}
\left\{
\begin{aligned}
\psi_i^l & = \arg\min_{\m} F[\m] \\
& = \arg\min_{\m} -\frac{1}{2} \int_{\Omega_i^l} \h_{\eff}[\m] \cdot \m \d\x \\
& = \arg\min_{\m} \frac{1}{2} \int_{\Omega_i^l} \kappa |\nabla \m|^2 + m_2^2 + m_3^2 \d\x, \\
& \text{s.t.} \quad \langle \phi_j, \m \rangle = \delta_{ij}, \quad \forall \phi_j \in \Phi^l,
\end{aligned}
\right.
\end{equation}
where $\Phi^l$ denotes the measurement functions in local patch $\Omega_i^l$.

\begin{remark}
For the local optimization problem \eqref{vector_optimize} on $\Omega_i^l \subsetneq \Omega$, homogeneous Dirichlet conditions hold on $\partial\Omega_i^l$ when $\Omega_i^l \subset \Omega$. For subdomains intersecting the global boundary ($\partial\Omega_i^l \cap \partial\Omega \neq \emptyset$), mixed conditions apply with Neumann on $\partial\Omega_i^l \cap \partial\Omega$ and Dirichlet on $\partial\Omega_i^l \setminus \partial\Omega$. 
\end{remark}

Given the vectorial nature of the equation and the anisotropy alignment $\u=(1,0,0)^T$, we solve two different variational problems from \eqref{vector_optimize}:

\begin{equation*}\label{basis_cp1}
\left\{
\begin{aligned}
\psi_i^l &= \arg\min_{v} \int_{\Omega_i^l} \kappa |\nabla v|^2 \d\x \\
&\text{s.t.} \quad \langle \phi_j,v \rangle = \delta_{ij}, \quad \forall \phi_j \in \Phi^l,
\end{aligned}
\right.
\tag{V1}
\end{equation*}

\begin{equation*}\label{basis_cp23}
\left\{
\begin{aligned}
\psi_i^l &= \arg\min_{v} \int_{\Omega_i^l} (\kappa |\nabla v|^2 + v^2) \d\x \\
&\text{s.t.} \quad \langle \phi_j,v \rangle = \delta_{ij}, \quad \forall \phi_j \in \Phi^l.
\end{aligned}
\right.
\tag{V2}
\end{equation*}

The basis from (V1) spans the GRPS space for $m_1$, while (V2) generates the basis for $m_2$ and $m_3$ components, reflecting the magnetic-anisotropic energy structure.

\begin{remark}
Through the utilization of the energy minimization approach for the LL equation, there is no need to compute the multiscale basis at each time step. Besides, as documented in references \cite{zhang_RPS,liuGRPS}, an additional advantage of this method lies in the parallel computation of the multiscale basis, leading to significant savings in computational cost and time.
\end{remark}

\paragraph{Exponential decay of multiscale basis functions}

We demonstrate the exponential decay properties of the multiscale basis functions through numerical experiments. The spatial domain $\Omega = [0, 1]^2$ is partitioned into $N_c = 2^3$ patches, each subdivided into $J = 4$ segments, yielding a fine mesh size $h = 1/2^7$.

\begin{figure}[H]
    \centering
    \subfigure[\eqref{basis_cp1}, GRPS-E]{
        \includegraphics[width=4.4cm]{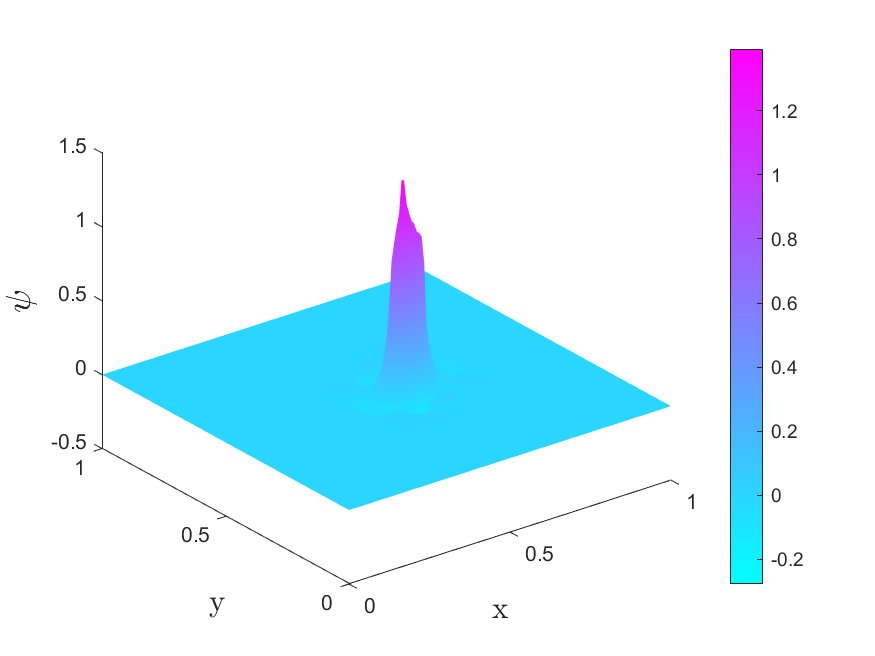}
        \label{Fig1-1}
    }
    \quad
    \subfigure[\eqref{basis_cp1}, GRPS-V]{
        \includegraphics[width=4.4cm]{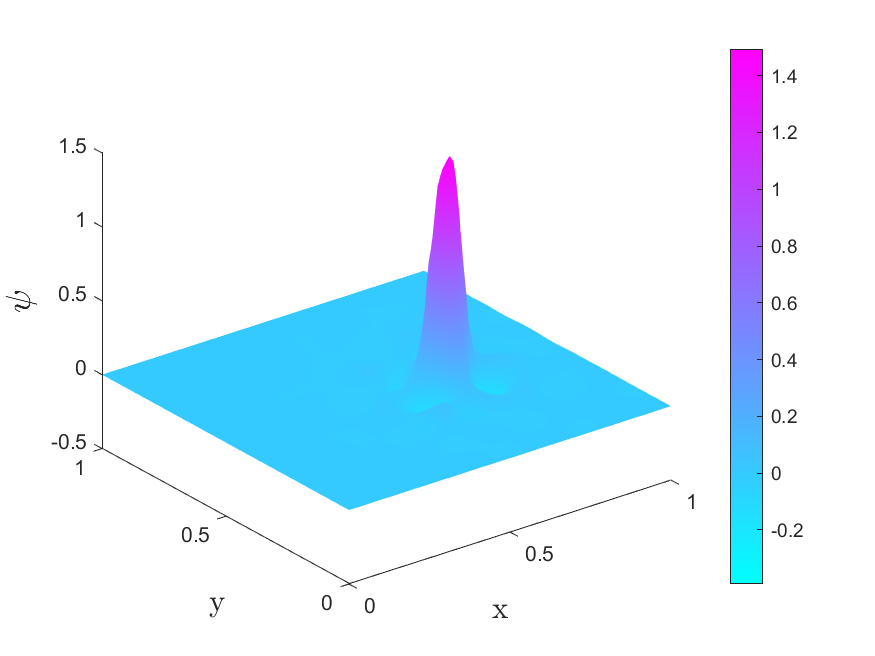}
        \label{Fig1-2}
    }
    \quad    
    \subfigure[\eqref{basis_cp23}, GRPS-E]{
        \includegraphics[width=4.4cm]{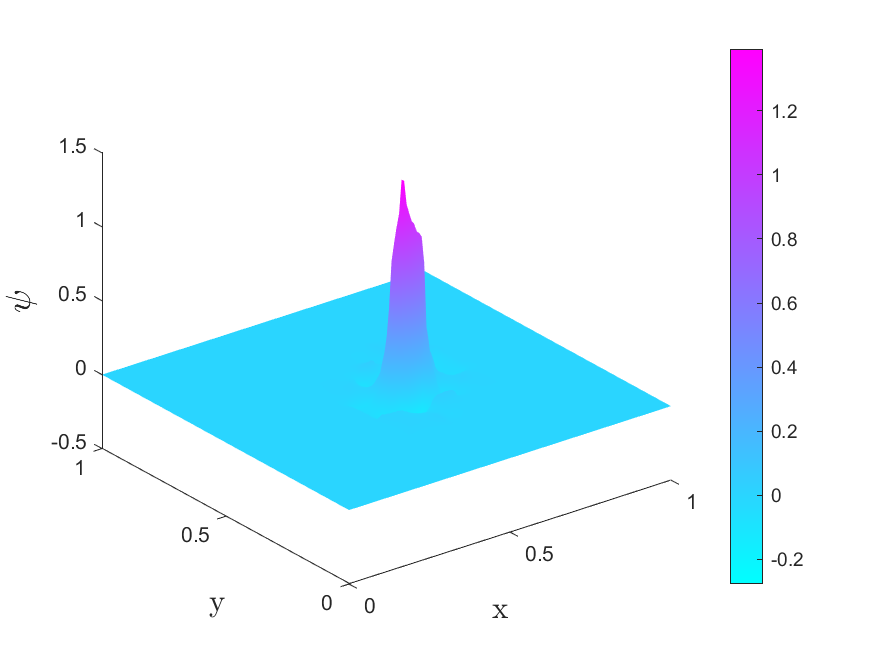}
        \label{Fig2-1}
    }
    \quad
    \subfigure[\eqref{basis_cp23}, GRPS-V]{
        \includegraphics[width=4.4cm]{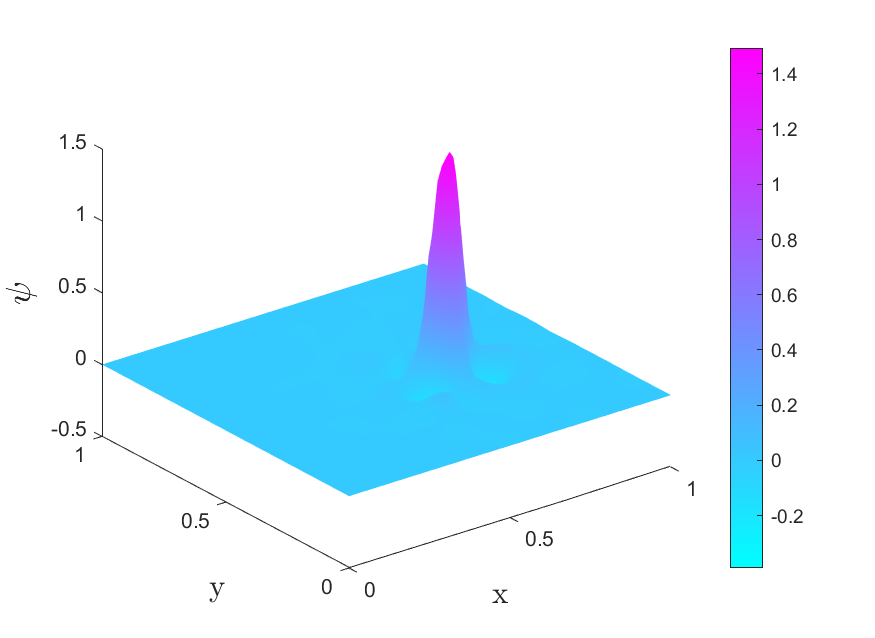}
        \label{Fig2-2}
    }
    \caption{GRPS-E and GRPS-V basis functions obtained from \eqref{basis_cp1} and \eqref{basis_cp23}.} 
    \label{fig:profileV1andV2}
\end{figure}

\begin{figure}[H]
    \centering
    \subfigure[\eqref{basis_cp1}, GRPS-E]{
        \includegraphics[width=4.4cm]{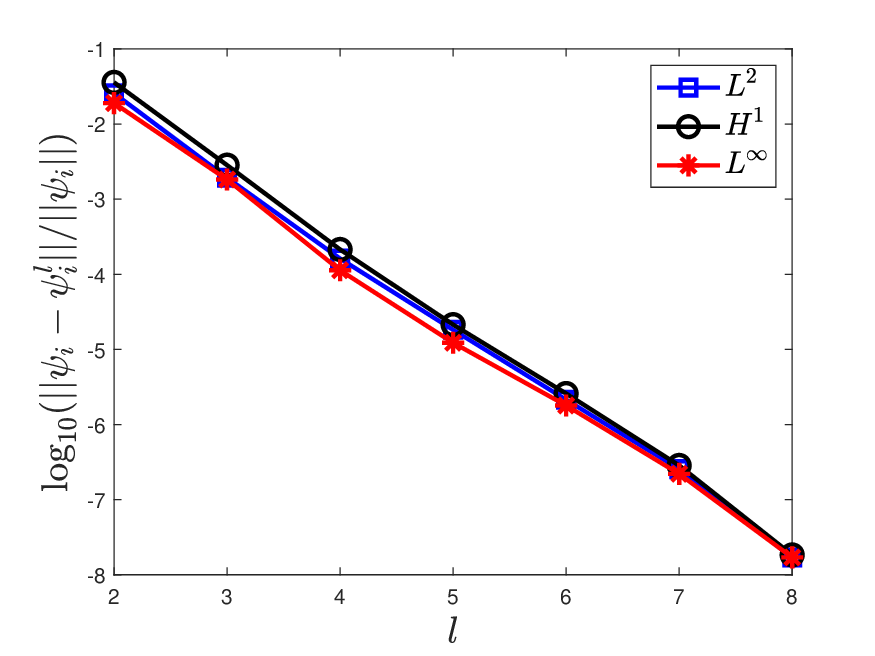}
        \label{Fig1-1expdecay}
    }
    \quad
    \subfigure[\eqref{basis_cp1}, GRPS-V]{
        \includegraphics[width=4.4cm]{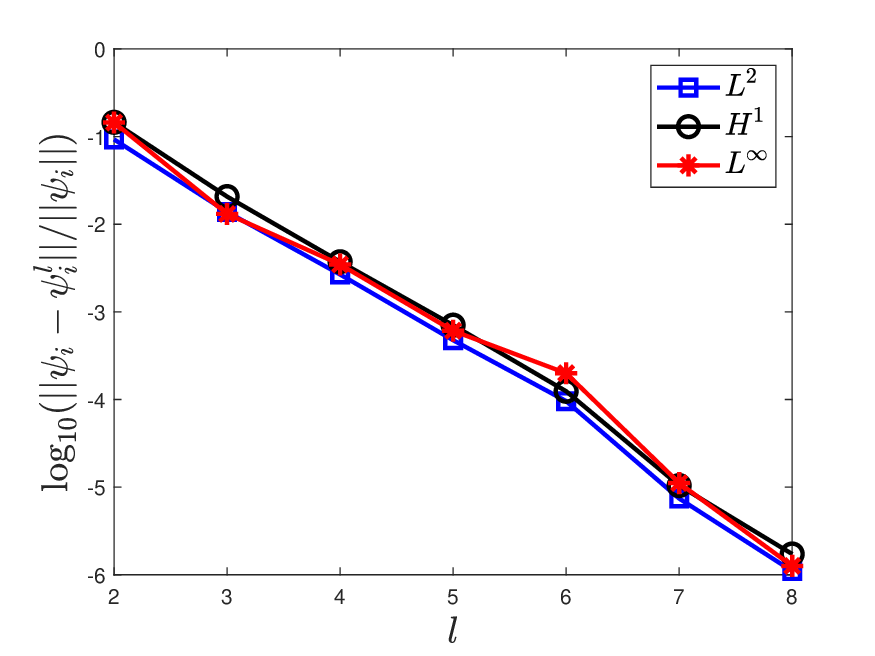}
        \label{Fig1-2expdecay}
    }
    \quad    
    \subfigure[\eqref{basis_cp23}, GRPS-E]{
        \includegraphics[width=4.4cm]{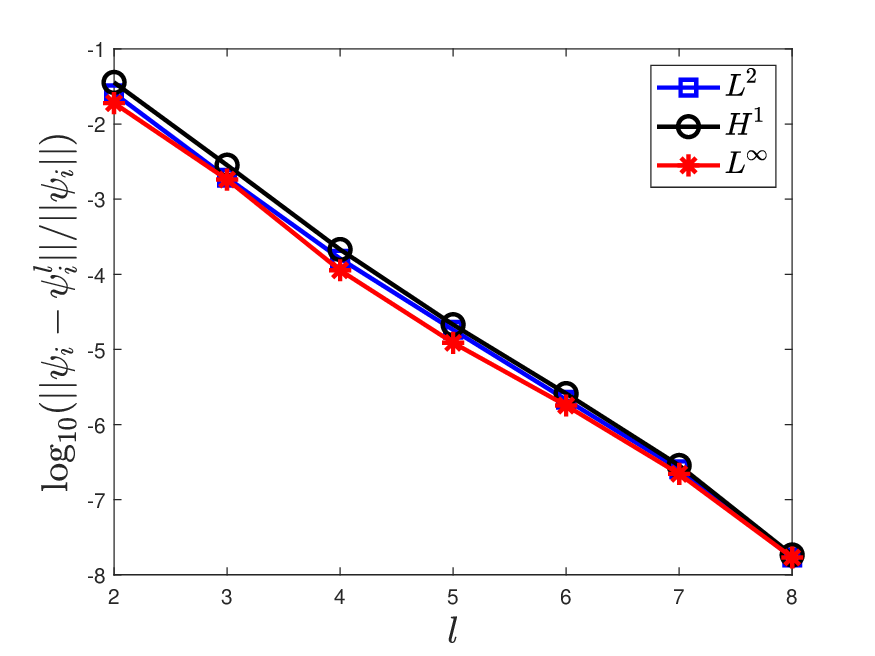}
        \label{Fig2-1expdecay}
    }
    \quad
    \subfigure[\eqref{basis_cp23}, GRPS-V]{
        \includegraphics[width=4.4cm]{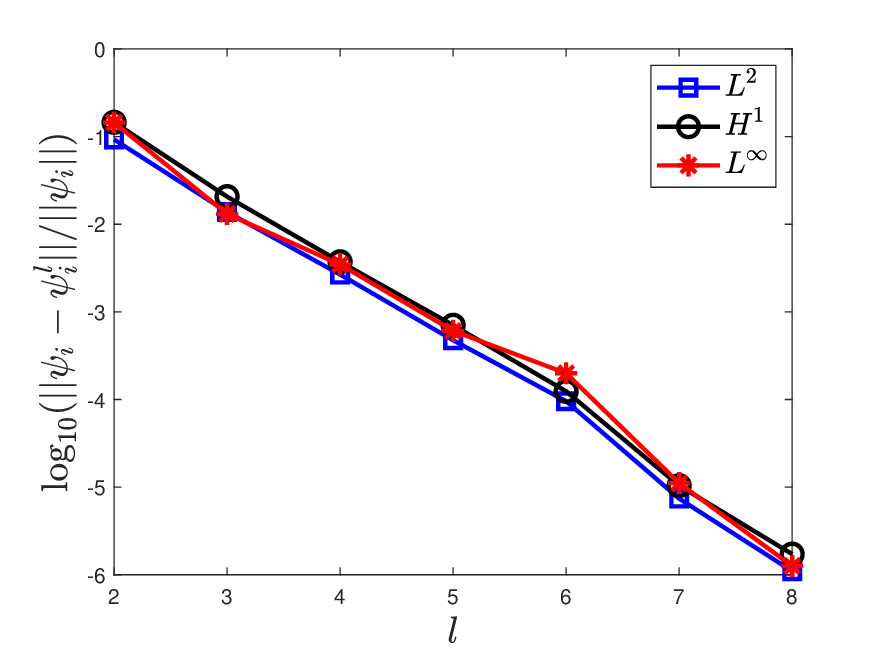}
        \label{Fig2-2expdecay}
    }
    \caption{Exponential decay measured by $||\psi_i - \psi_i^l||/||\psi_i||$ in $L^2$, $H^1$, and $L^{\infty}$ norms for localization levels $l=2,\cdots,8$ (log$_{10}$ scale).}
    \label{fig:expdecay}
\end{figure}

Figures \ref{fig:profileV1andV2} and \ref{fig:expdecay} display the profiles and decay characteristics of GRPS-E and GRPS-V basis functions derived from \eqref{basis_cp1} and \eqref{basis_cp23}, with $\kappa$ given by \eqref{multiscale_coefficient_k}. The exponential decay is quantified by plotting $\log(||\psi_i - \psi_i^l||/||\psi_i||)$ against localization level $l$, confirming the expected decay behavior for both basis types.

\subsection{Multiscale Algorithms for Landau-Lifshitz Equation}

With the coarse approximation space established, we now develop two distinct numerical approaches for solving the LL equation \eqref{LLform1}. The magnetization field $\m$ in micromagnetics must theoretically satisfy the constraint $|\m| = 1$ at all points. While our first algorithm \ref{alg:length relaxing} relaxes this constraint, we demonstrate through numerical analysis \cite{prohl2001computational,cimrak2005error,gao2014optimal} that careful discretization can maintain $|\m| \approx 1$ within controlled error bounds even without explicit enforcement. A \textit{strict-constraint} approach (Algorithm \ref{alg:length preserving}) explicitly enforces $|\m| = 1$ at each computational step.

\begin{algorithm}[H]
\begin{algorithmic}[1]
\STATE Obtain the multiscale basis $\Psi = \{ \psi_j \}_{j=1}^{N_H}$ by solving \eqref{vector_optimize}.
\STATE Prescribe initial conditions: $\m_h^0$, $\m_H^0= \sum \limits_{j} \langle \phi_{j}, \m_h^0 \rangle \psi_j$.
\FOR{$n=0$ to $N-1$}
\STATE Obtain coarse scale solution $\m_{H}^{n+1}$ by solving \eqref{nospeed1}.
\ENDFOR
\end{algorithmic}
\caption{Multiscale length relaxing backward FEM for LL equation}
\label{alg:length relaxing}
\end{algorithm}

\begin{algorithm}[H]
\begin{algorithmic}[1]
\STATE Obtain the multiscale basis $\Psi = \{ \psi_j \}_{j=1}^{N_H}$ by solving \eqref{vector_optimize}.
\STATE Prescribe initial conditions $\m_h^0$, $\tilde{\m}_H^0= \sum \limits_{j} \langle \phi_{j}, \m_h^0 \rangle \psi_j$ .
\FOR{$n=0$ to $N-1$}
\STATE Projection onto unit sphere $\m_{H}^{n}=\frac{ \tilde{\m}_{H}^{n}}{| \tilde{\m}_{H}^{n}|}$;
\STATE Obtain the intermediate coarse scale solution $\tilde{\m}_{H}^{n+1}$ by solving \eqref{nospeed1}.
\ENDFOR
\STATE Projection onto unit sphere $\m_{H}^{N}=\frac{ \tilde{\m}_{H}^{N}}{| \tilde{\m}_{H}^{N}|}$.
\end{algorithmic}
\caption{Multiscale length preserving backward FEM for LL equation}
\label{alg:length preserving}
\end{algorithm}
In Algorithm \ref{alg:length preserving}, it is important to note that the preservation of unit length on the fine scale is implied by $\m_{H}^{n}=\frac{\tilde{\m}_{H}^{n}}{| \tilde{\m}_{H}^{n}|}$ rather than projection step of coefficients on the coarse scale directly, with $ \tilde{\m}_{H}^{n}$ indicates the fine scale intermediate solution. This choice is based on the fact that $\m_{H}^{n}$ provides a good approximation of $\m_h^n$ at the fine scale, while the direct division operation at the coarse scale eliminates the scaling factor $c_e$ or $c_{\tau}$, which will introduce issues in the subsequent time step evolution.





\subsection{Speed-up Operation for Nonlinear Term Assembly}
The assembly of nonlinear terms using classical FEM is a time-intensive process. For instance, the 4-valence tensor necessitates four iterations to traverse each fine simplex, as illustrated by the following equation:
\begin{equation*}
    (|\nabla \m_h^{j}|^2 \varphi_j, \varphi_i ) = ((\sum \limits _{k=1}^{N_h} \m_{h,k}^{j} \nabla \varphi_k)  (\sum \limits _{l=1}^{N_h} \m_{h,l}^{j} \nabla \varphi_l) \varphi_j, \varphi_i ),
\end{equation*}
where $\varphi_i, \varphi_j, \varphi_k, \varphi_l \in V_h$ with $i,j,k,l = 1, \cdots, N_h.$ So as the 4-valence tensor for the anisotropic term $((\sum \limits _{k=1}^{N_h} \m_{h,k}^{j} \varphi_k)  (\sum \limits _{l=1}^{N_h} \m_{h,l}^{j} \varphi_l) \varphi_j, \varphi_i )$. Besides, when dealing with a fourth-order polynomial, a minimum of a 6-point quadrature formula becomes essential. 

The backward Euler schemes \eqref{eqn:fullynonlinear}-\eqref{AnFEM}, particularly \eqref{Cimark} and \eqref{GaoFEM}, enable computational acceleration through their square term incorporation. However, pre-computing the required 4-valence tensor in GRPS space exceeds practical memory limits due to dimensionality, which motivates us to transform the 4-valence tensor into a 3-valence tensor following \cite{henning2022superconvergence}.

For scheme \eqref{Cimark} with GRPS bases, we seek $\m_{H}^{n+1} \in \V_H$ satisfying
\begin{equation}\label{nospeed1}
    A^n(\m_{H}^{n+1}, \v) = (\f_H^n,\v_H), \quad \forall \v_H \in \V_H.
\end{equation}

Define the $L^2$-projection $P_{GRPS}: H^1(\Omega) \rightarrow V_H$,
\begin{equation}
    (P_{GRPS}(u),v) = (u,v), \quad \forall v \in V_{H}.
\end{equation}

we propose an accelerated variant,
\begin{equation}\label{speed}
    \begin{aligned}
        &(D_{\tau} \m_{H}^{n+1},\v) + B^n(\m_H^{n+1},\v) - (\m_{H}^{n} \times \m_{H,a}^{n+1}, \v) \\
        =& \lambda (\kappa(\x) P_{GRPS}(|\nabla \m_{H}^{n}|^2 ) \m_{H}^{n+1}, \v) + \lambda (P_{GRPS}( |m_{H,2}^{n}|^2 + |m_{H,3}^{n}|^2 ) \m_{H}^{n+1},\v).
    \end{aligned}
\end{equation}

The acceleration stems from reduced assembly complexity,
\begin{equation}
    (\kappa(\x) P_{GRPS}(|\nabla \m_{H}^{n}|^2 ) \m_{H}^{n+1},\psi_l) = (\kappa(\x) \sum_{i,j=1}^{N_H} \rh_{i}^{n} \m_{H,j}^{n+1} \psi_i \psi_j, \psi_l).
\end{equation}

Define tensor $\boldsymbol{\omega} = (\omega_{kji})$ and vector $\dd = (d_i)$,
\begin{equation}
    \omega_{kji} = (\partial_x \psi_k \partial_x \psi_j + \partial_y \psi_k \partial_y \psi_j, \psi_i),
\end{equation}
\begin{equation}
    d_i = (|\nabla \m_{H}^{n}|^2, \psi_i).
\end{equation}

The coefficients $\rh_{l}^{n}$ solve the following linear system,
\begin{equation}
    \dd= M \rh^{n}, \quad M_{il} = (\psi_l,\psi_i).
\end{equation}

The final assembly uses precomputed tensors $\boldsymbol{\omega}$ and $\boldsymbol{\Bar{\omega}}$, both facilitating reusable computations, with analogous treatment for the anisotropy term $\boldsymbol{\Bar{\omega}}$:
\begin{equation}
    \Bar{\omega}_{kji} = ( \psi_k \psi_j, \psi_i).
\end{equation}

\section{Numerical Experiments}
\label{sec4}

We conduct numerical tests to evaluate our multiscale method's performance, specifically examining the variational formulation setup, treatment of the length-preserving step, and acceleration of nonlinear terms.

\paragraph{Experimental Setup}
All 2D simulations use $\Omega = [0,1]^2$ with $T=1.0$. Initial conditions are $\m_0 = \left(1/\sqrt{2}, 1/\sqrt{3}, 1/\sqrt{6}\right)^T$ with homogeneous Neumann boundary conditions and damping parameter $\lambda=1$.

We measure solution accuracy using relative errors at final time $T$:
\begin{equation}
\frac{\|\m_h^{N} - \m_H^{l,N}\|_{\cdot}}{\|\m_h^{N}\|_{\cdot}},
\end{equation}
where $\|\cdot\|$ denotes $H^1$ norm, $l$ is the localization level, and $N$ represents the final time step.

\begin{figure}[htbp]
\centering
\includegraphics[height=5cm,width=6cm]{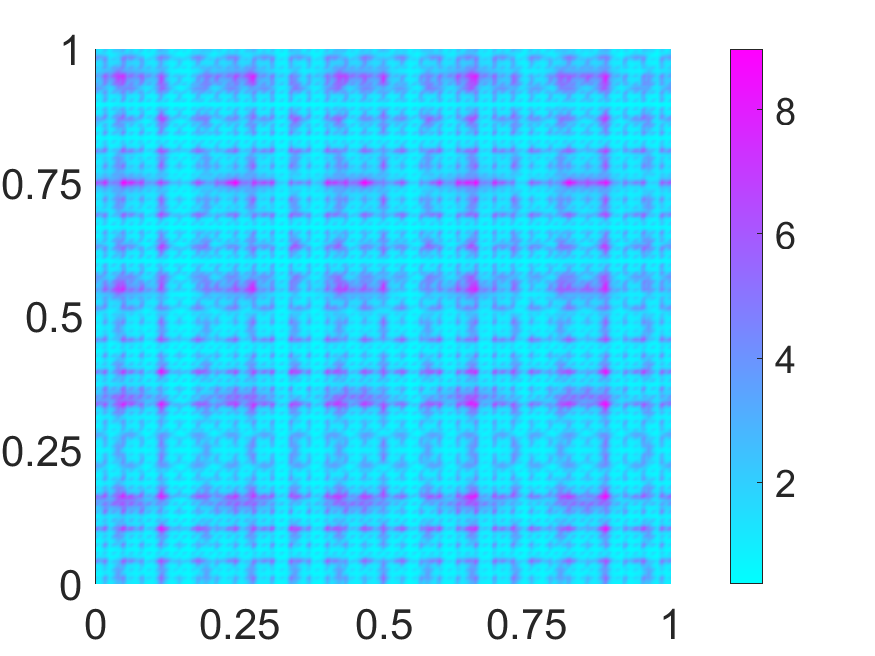}
\caption{MsTrig coefficient $\kappa(x,y)$ on $\Omega=[0,1]^2$.}
\label{fig_kappa}
\end{figure}

\paragraph{Multiscale Coefficient}
We employ the oscillatory MsTrig coefficient:
\begin{equation}\label{multiscale_coefficient_k}
\kappa(x,y) = \frac{1}{6}\left[
\sum_{i=1}^5 f_i(x,y,\eps_i) + \sin(4x^2y^2) + 1
\right],
\end{equation}
where $\eps_i = \{1/5, 1/13, 1/17, 1/31, 1/65\}$ and $f_i$ are trigonometric terms (see Appendix \ref{subsec:six-scale multiscale coefficient}). Figure \ref{fig_kappa} shows its highly oscillatory profile.

Following \cite{zhang_RPS,owhadi2017multigrid,liuGRPS}, we expect error bounds consistent with Theorem \ref{theorem_convergence_rate} when applying GRPS bases to the LL equation.

\paragraph{Key Parameters}
Key factors under examination:
\begin{itemize}
    \item Basis variational forms (V1, V2);
    \item Measurement functions (E, V);
    \item Coarse time step choices ($\tau_H = H$ vs $H^2$);
    \item Length-preservation implementation.
\end{itemize}

\subsection{Accuracy Tests: The Influence of the Key Parameters}

We employ the schemes described in Appendix \ref{subsec:Various Time Discretization Schemes}, as the length-preserving step has been rigorously justified theoretically. The non-separable multiscale coefficient \eqref{multiscale_coefficient_k} is used with final time $T$, spatial domain $\Omega$, and damping parameter $\lambda$ ($T$, $\Omega$, $\lambda$) consistent with previous settings. The reference solution on a fine mesh uses $h=1/2^{7}$ and $\Delta t=h^2$, while coarse grids are divided into $N_c = 2, 4, 8, 16$ ($H=1/N_c$) with time steps $\tau_H = H, H^2$.

\subsubsection{Without Length Preserving Step}
The effective field $\h_{\eff} = \div(\kappa \nabla \m) - (\m - (\m \cdot \u)\u)$ remains crucial, as established by Landau and Lifshitz \cite{landau1992theory}. The multiscale basis from optimization problems \eqref{basis_cp1}, \eqref{basis_cp23}, and their combination maintains effectiveness with the anisotropy term.

In Tables \ref{table_with_anisotrpy}-\ref{Gao_table_without_lengthpre}, we record the global relative $H^1$ error convergence rates on a series of coarse degrees of freedoms when $N_c = 2, 4, 8, 16$ with diverse types of measurement functions, basis variational forms and coarse time steps. Upon examination of An's method, Cimr\'ak's method, and Gao's method, it is observed that for coarse time steps $\tau_H = H, H^2$, both the proposed GRPS-E and GRPS-V methods converge with rates of 1 and 2, respectively. Given that the previous three methods exhibit first-order accuracy in time and at least first-order accuracy in space within the multiscale Finite Element Method (FEM) framework concerning $H^1$ norm (with $O(H^2)$ accuracy for GRPS-V), it is evident that the selection of coarse time steps could impact the spatial order. However, for GRPS-E, which exhibits higher order accuracy, this may be attributed to the super convergence phenomenon in this problem.

\begin{table*}[htbp]
    \centering
    \caption{$\kappa = MsTrig$. An's method. Global convergence rates comparing fine mesh reference solution to multiscale solution using GRPS-E and GRPS-V bases. TMB stands for "Type of multiscale bases," while TMF stands for "Type of measurement function." Fine mesh: $h=1/2^{7}$, $\Delta t=h^2$. Coarse grid: $N_c =2,4,8,16$, $\tau_H = H, H^2$.}
    \label{table_with_anisotrpy}
    \begin{tabular}{c|ccccc}
        \Xhline{2pt}
         \multirow{2}{*}{\diagbox{TMB}{TMF}} & \multicolumn{4}{c}{Convergence Rate} \\
         & GRPS-E,$\tau_H=H$ & GRPS-E,$\tau_H=H^2$ & GRPS-V,$\tau_H=H$ & GRPS-V,$\tau_H=H^2$  \\
        \Xhline{1pt}
        \eqref{basis_cp1}+\eqref{basis_cp23} & -0.4182 & -0.8187 & -0.4643 & -0.9940 \\
        \eqref{basis_cp1} & -0.4303 & -0.9174 & -0.4833 & -0.9940 \\  
        \eqref{basis_cp23} & -0.4172 & -0.6646 & -0.4833 & -0.9941 \\  
        \Xhline{2pt}
    \end{tabular} 
\end{table*}

\begin{table*}[htbp]
    \centering
    \caption{Cimr\'ak's method. Global convergence rates without length preserving constraint. Parameters as in Table \ref{table_with_anisotrpy}.}
    \label{Cimrak_table_without_lengthpre}
    \begin{tabular}{c|ccccc}
        \Xhline{2pt}
         \multirow{2}{*}{\diagbox{TMB}{TMF}} & \multicolumn{4}{c}{Convergence Rate} \\
         & GRPS-E,$\tau_H=H$ & GRPS-E,$\tau_H=H^2$ & GRPS-V,$\tau_H=H$ & GRPS-V,$\tau_H=H^2$  \\
        \Xhline{1pt}
        \eqref{basis_cp1}+\eqref{basis_cp23} & -0.51730 & -0.89267 & -0.55234 & -1.0325 \\
        \eqref{basis_cp1} & -0.51422 & -0.95192  & -0.52973 & -1.0325 \\  
        \eqref{basis_cp23} & -0.51675 & -0.74907 & -0.52672 & -1.0210 \\  
        \Xhline{2pt}
    \end{tabular} 
\end{table*}

\begin{table*}[htbp]
    \centering
    \caption{Gao's method. Global convergence rates without length preserving constraint. Parameters as in Table \ref{table_with_anisotrpy}.}
    \label{Gao_table_without_lengthpre}
    \begin{tabular}{c|ccccc}
        \Xhline{2pt}
         \multirow{2}{*}{\diagbox{TMB}{TMF}} & \multicolumn{4}{c}{Convergence Rate} \\
         & GRPS-E,$\tau_H=H$ & GRPS-E,$\tau_H=H^2$ & GRPS-V,$\tau_H=H$ & GRPS-V,$\tau_H=H^2$  \\
        \Xhline{1pt}
        \eqref{basis_cp1}+\eqref{basis_cp23} & -0.44480 & -0.86602 & -0.48332 & -0.99940 \\
        \eqref{basis_cp1} & -0.44427 & -0.92084 & -0.48332 & -0.99940 \\  
        \eqref{basis_cp23} & -0.44847 & -0.73035 & -0.48332 & -0.99940 \\  
        \Xhline{2pt}
    \end{tabular} 
\end{table*}

\begin{figure}[htbp]
\centering
\subfigure[GRPS-E: \eqref{basis_cp1}+\eqref{basis_cp23}]{
\includegraphics[scale=0.45]{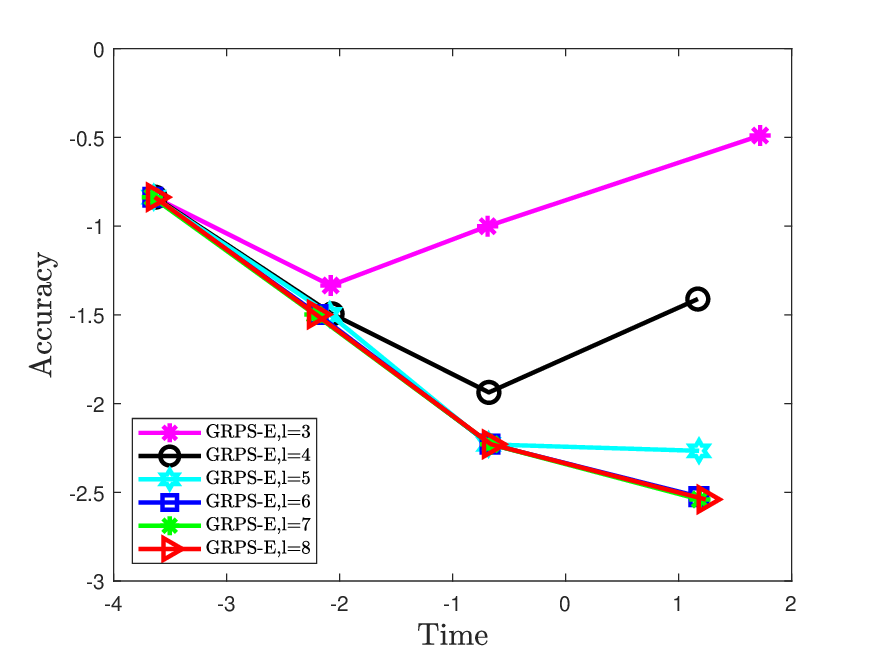}
}
\hfill
\subfigure[GRPS-V: \eqref{basis_cp1}+\eqref{basis_cp23}]{
\includegraphics[scale=0.45]{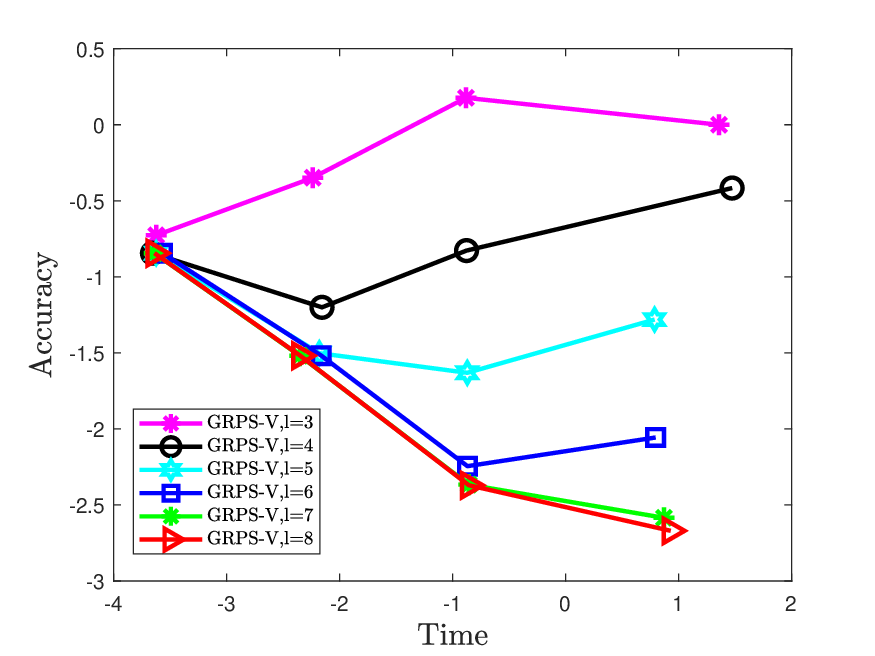}
}
\caption{Performance of GRPS-E and GRPS-V multiscale bases solving \eqref{basis_cp1}+\eqref{basis_cp23} with localization levels $l=3,4,5,6,7,8$. X-axis: time (seconds), Y-axis: $H^1$ relative error (log10 scale).}
\label{fig:time_accuracy_withaniso_caseET}
\end{figure}

Next, we compare the time-accuracy of GRPS-E and GRPS-V on a laptop with an Intel i9-14900HX CPU at 2.2GHz $\times$ 32 processors and 64GB of RAM running MATLAB version R2023a, as shown in Figure \ref{fig:time_accuracy_withaniso_caseET}. Since the computation time is proportional to the degrees of freedom (dofs) of the coarse mesh, the time-accuracy table also exhibits the localization properties, similar to the dof-accuracy table, with the increasing number of the localization parameter $l$ from 3 to 8 in Figure \ref{fig:time_accuracy_withaniso_caseET}.

To achieve nearly the same accuracy on both fine and coarse meshes, i.e., $h=H^2$, we set $\tau_H = H^2$ on the coarse scale to achieve $O(H^2)$ accuracy in the $H^1$ norm for spatial accuracy. Correspondingly, on the fine mesh, we set $\Delta t = h$ to attain $O(h)$ accuracy. Numerical results suggest that the reference solution ($h=2^{-7}$, $\Delta t=h$) required 36 seconds (excluding the assembly of the 4-valence tensors, which is extremely time-consuming in nonlinear problems, focusing solely on solving the linear systems). In comparison, GRPS-E ($N_c = 16$) took 15.4 seconds (a reduction of 57.22\%), and GRPS-V took 7.4 seconds (a reduction of 79.44\%). Both GRPS-E and GRPS-V show significant efficiency improvements when employing the numerical scale-upscaling procedure. The main difference in computational time between GRPS-E and GRPS-V arises from the increased number of coarse degrees of freedom in GRPS-E compared to GRPS-V, resulting in higher computational costs.

\subsubsection{With Length Preserving Step}
Incorporating the unit length constraint (Algorithm \ref{alg:length preserving}) and exchange/anisotropy terms, we use coefficient \eqref{multiscale_coefficient_k} with $h=1/2^{7}$, $\Delta t=h^2$, and $N_c=2,4,8,16$ ($\tau_H=H,H^2$). 

\begin{table*}[htbp]
    \centering
    \caption{An's method. Global convergence rates with length preserving constraint. Parameters as in Table \ref{table_with_anisotrpy}.}
    \label{An_table_with_anisotrpy_lengthpre}
    \begin{tabular}{c|ccccc}
        \Xhline{2pt}
         \multirow{2}{*}{\diagbox{TMB}{TMF}} & \multicolumn{4}{c}{Convergence Rate} \\
         & GRPS-E,$\tau_H=H$ & GRPS-E,$\tau_H=H^2$ & GRPS-V,$\tau_H=H$ & GRPS-V,$\tau_H=H^2$ \\
        \Xhline{1pt}
        \eqref{basis_cp1}+\eqref{basis_cp23} & -0.49115 & -0.79407 & -0.50260 & -1.0028 \\
        \eqref{basis_cp1} & -0.46229 & -0.92290 & -0.50260 & -1.0028 \\  
        \eqref{basis_cp23} & -0.46230 & -0.92312 & -0.50260 & -1.0028 \\  
        \Xhline{2pt}
    \end{tabular}
\end{table*}

\begin{table*}[htbp]
    \centering
    \caption{Cimr\'ak's method. Global convergence rates with length preserving constraint. Parameters as in Table \ref{table_with_anisotrpy}.}
    \label{cimrak_table_with_anisotrpy_lengthpre}
    \begin{tabular}{c|ccccc}
        \Xhline{2pt}
         \multirow{2}{*}{\diagbox{TMB}{TMF}} & \multicolumn{4}{c}{Convergence Rate} \\
         & GRPS-E,$\tau_H=H$ & GRPS-E,$\tau_H=H^2$ & GRPS-V,$\tau_H=H$ & GRPS-V,$\tau_H=H^2$ \\
        \Xhline{1pt}
        \eqref{basis_cp1}+\eqref{basis_cp23} & -0.46975 & -0.72577 & -0.50895 & -0.99859 \\
        \eqref{basis_cp1} & -0.46830 & -0.91898 & -0.50895 & -0.99859 \\  
        \eqref{basis_cp23} & -0.46614 & -0.91770 & -0.50895 & -0.99859 \\  
        \Xhline{2pt}
    \end{tabular}
\end{table*}

\begin{table*}[htbp]
    \centering
    \caption{Gao's method. Global convergence rates with length preserving constraint. Parameters as in Table \ref{table_with_anisotrpy}.}
    \label{Gao_table_with_anisotrpy_lengthpre}
    \begin{tabular}{c|ccccc}
        \Xhline{2pt}
         \multirow{2}{*}{\diagbox{TMB}{TMF}} & \multicolumn{4}{c}{Convergence Rate} \\
         & GRPS-E,$\tau_H=H$ & GRPS-E,$\tau_H=H^2$ & GRPS-V,$\tau_H=H$ & GRPS-V,$\tau_H=H^2$ \\
        \Xhline{1pt}
        \eqref{basis_cp1}+\eqref{basis_cp23} & -0.46185 & -0.79233 & -0.50260 & -1.0028 \\
        \eqref{basis_cp1} & -0.46229 & -0.92290 & -0.50260 & -1.0028 \\  
        \eqref{basis_cp23} & -0.46230 & -0.92313 & -0.50260 & -1.0028 \\  
        \Xhline{2pt}
    \end{tabular}
\end{table*}

In Tables \ref{An_table_with_anisotrpy_lengthpre}-\ref{Gao_table_with_anisotrpy_lengthpre}, we employ length preserving Algorithm \ref{alg:length preserving}, we record the global relative $H^1$ error convergence rates on a range of coarse degrees of freedom when $N_c = 2, 4, 8, 16$ with diverse types of measurement functions, basis variational forms and coarse time steps. Similar to what was observed in Tables \ref{table_with_anisotrpy}-\ref{Gao_table_without_lengthpre}, Tables \ref{An_table_with_anisotrpy_lengthpre}-\ref{Gao_table_with_anisotrpy_lengthpre} illustrate that the proposed bases are all effective, and the global convergence rates remain consistent with those presented above, even without length preserving considerations.


\paragraph{Summary}
After solving the LL equation using various backward Euler schemes—including An's method, Cimrák's method, and Gao's method — we observe that all methods achieve the same order of accuracy and exhibit nearly identical computational efficiency when applied with different basis variational forms (Eqs. \eqref{basis_cp1} and \eqref{basis_cp23}). Additionally, as demonstrated in prior tests, the selection of coarse time steps can negatively impact the spatial error order.

Furthermore, proper enforcement of the length-preservation step does not significantly affect runtime but is crucial for maintaining the unit-length constraint. Regarding measurement functions, a comparison between GRPS-E and GRPS-V reveals that GRPS-V is computationally more efficient, as it requires fewer basis functions than GRPS-E, leading to reduced runtime.

\subsection{The Performance of Speed-Up Operation}
Since GRPS-V is more efficient than GRPS-E in terms of computational time, we compare the computational time between classical FEM and GRPS-V basis for equivalent $H^1$ error orders, using identical parameters from previous tests. Table \ref{GRPS_T_timecompare2} shows:

\begin{itemize}
    \item CT0: Reference solution time (finest mesh);
    \item CT1: Fine matrix assembly (4-valence tensor, days) + CT0;
    \item CT2: Multiscale basis generation;
    \item CT3: 3-valence tensor assembly; 
    \item CT4: Coarse scale solution time.
\end{itemize}

For coarse time step $\tau_H = H^2$, GRPS-V achieves second-order $H^1$ accuracy in space. At $N_c=2^4$ ($h = H^2$), the multiscale approach matches FEM precision $O(h)$, resulting in a significant reduction in time from CT0 to CT4, greatly reducing runtime by $\geq$94.4\%.

\begin{table}[H]
    \centering
    \caption{Time comparison between standard \eqref{nospeed1} and accelerated \eqref{speed} schemes in GRPS-V space. Parameters: $h=1/2^8$, $\Delta t=h$; coarse grid $N_c=2,4,8,16$ with $\tau_H=H^2$.}
    \label{GRPS_T_timecompare2}
    \begin{tabular}{c|cccc}
        \Xhline{2pt}
        \diagbox{Time}{$N_c$} & 2 & 4 & 8 & 16 \\
        \Xhline{1pt}
        CT0 (s) & 418 & 418 & 418 & 418 \\
        CT1 (Day) & 3.1 & 3.1 & 3.1 & 3.1 \\
        \Xhline{1pt}
        CT2 (s) & 0.176 & 1.225 & 9.414 & 27 \\
        CT3 (s) & 673 & 3589 & 10573 & 14951 \\
        CT4 (s) & 2.04e-4 & 6.81e-3 & 1.32e-1 & 6.99 \\
        \Xhline{2pt}
    \end{tabular}
\end{table}
Table \ref{GRPS_T_timecompare2} illustrates that the assembly of the finest 4-valence tensor (CT1) is significantly time and memory-intensive due to the curse of dimensionality. Implementing speed-up operations notably decreases both assembly time and memory requirements (CT3). CT2 indicates the generation of the GRPS-V bases, each of which is independent and can be solved concurrently. Furthermore, utilizing GRPS-V bases facilitates solving linear systems in a low-dimensional space, resulting in a significant reduction in time from CT0 to CT4. Finally, the speed-up operations achieve nearly the same level of accuracy as shown in Table \ref{cimrak_table_with_anisotrpy_lengthpre}. Therefore, a detailed display will not be provided here.
\begin{remark}
    The accuracy of the accelerated scheme \eqref{speed} is almost the same with the Tables \ref{table_with_anisotrpy}-\ref{Gao_table_without_lengthpre} under examinations, and some related rigorous proofs are detailed in \cite{henning2022superconvergence}. Hence, we will not include an excessive number of tables and instead focus solely on the performance of accelerated schemes.
\end{remark}

\section{Conclusion and Discussion}
\label{sec:conclusion}

This paper presents a numerical homogenization method for solving the Landau-Lifshitz equation with rough coefficients. Our results demonstrate that when $\h_{\text{eff}}$ includes an anisotropy term, both GRPS-E and GRPS-V multiscale bases obtained through minimization problems \eqref{basis_cp1} and \eqref{basis_cp23} effectively capture the essential information of the divergence-form operator. 

From a physical perspective, the inclusion of anisotropy terms proves crucial for proper information extraction, consistent with the Landau-Lifshitz theory \cite{landau1992theory}. The model selection for the optimization problem significantly impacts our ability to compress operator information during multiscale basis construction. Computationally, our method achieves substantial efficiency gains compared to direct fine-mesh FEM solutions, while maintaining solution accuracy through coarse-mesh upscaling.

Future research directions will focus on three key aspects: First, establishing rigorous convergence analysis for the proposed method. Second, developing efficient numerical schemes for time-dependent coefficient cases. Third, extending the framework to handle more complex material configurations and boundary conditions.

\section*{Acknowledgments}

ZM, RD and LZ were partially supported by the National Natural Science Foundation of China (Grant No. 12271360). LZ was also partially supported by the Shanghai Municipal Science and Technology Project 23JC1402300 and the Fundamental Research Funds for the Central Universities. JC was partially supported by NSFC grant 12425113. The authors gratefully acknowledge Professor Rong An and Doctor Panchi Li for valuable discussions and insights.

\bibliographystyle{plain}
\bibliography{references}

@article{landau1992theory,
  author  = {L. D. Landau and E. M. Lifshitz},
  title   = {On the theory of the dispersion of magetic permeability in ferromagnetic bodies},
  journal = {Phys. Z. Sowjetunion}, 
  volume  = {8},
  year    = {1935},
  pages   = {153-169}
}

@article{gilbert1955lagrangian,
  title="{A Lagrangian formulation of the gyromagnetic equation of the magnetization field}",
  author={Gilbert, Thomas L},
  journal={Phys. Rev.},
  volume={100},
  pages={1243-1255},
  year={1955}
}

@article{hauck2023super,
  title={Super-localization of elliptic multiscale problems},
  author={Hauck, Moritz and Peterseim, Daniel},
  journal={Mathematics of Computation},
  volume={92},
  number={341},
  pages={981--1003},
  year={2023}
}

@book{prohl2001computational,
  title={Computational Micromagnetism},
  author={Prohl, Andreas and others},
  year={2001},
  publisher={Springer}
}

@article{cimrak2005error,
  title="{Error estimates for a semi-implicit numerical scheme solving the Landau--Lifshitz equation with an exchange field}",
  author={Cimr{\'a}k, Ivan},
  journal={IMA Journal of Numerical Analysis},
  volume={25},
  number={3},
  pages={611--634},
  year={2005},
  publisher={Oxford University Press}
}

@article{cimrak2007survey,
  title="{A survey on the numerics and computations for the Landau-Lifshitz equation of micromagnetism}",
  author={Cimr{\'a}k, Ivan},
  journal={Archives of Computational Methods in Engineering},
  volume={15},
  number={3},
  pages={1--37},
  year={2007},
  publisher={Springer}
}

@article{kruzik2006recent,
  title={Recent developments in the modeling, analysis, and numerics of ferromagnetism},
  author={Kruzik, Martin and Prohl, Andreas},
  journal={Siam Review},
  volume={48},
  number={3},
  pages={439--483},
  year={2006},
  publisher={SIAM}
}

@article{GarcaCervera2007NUMERICALMA,
  title={Numerical micromagnetics: a review},
  author={Cervera, Carlos J Garc{\'\i}a},
  journal={SeMA Journal: Bolet{\'\i}n de la Sociedad Espa{\~n}ola de Matem{\'a}tica Aplicada},
  volume={39},
  pages={103--135},
  year={2007},
  publisher={Sociedad Espa{\~n}ola de Matem{\'a}tica Aplicada}
}

@article{gao2014optimal,
  title="{Optimal error estimates of a linearized backward Euler FEM for the Landau--Lifshitz equation}",
  author={Gao, Huadong},
  journal={SIAM Journal on Numerical Analysis},
  volume={52},
  number={5},
  pages={2574--2593},
  year={2014},
  publisher={SIAM}
}

@article{parabolic_homo_zhang,
author = {Owhadi, Houman and Zhang, Lei},
title = {Homogenization of Parabolic Equations with a Continuum of Space and Time Scales},
journal = {SIAM Journal on Numerical Analysis},
volume = {46},
number = {1},
pages = {1-36},
year = {2008},
doi = {10.1137/060670420},
URL = { https://doi.org/10.1137/060670420}
}

@article{liuGRPS,
author = {Liu, Xinliang and Zhang , Lei and Zhu, Shengxin},
title = "{Generalized Rough Polyharmonic Splines for Multiscale PDEs with Rough Coefficients}",
journal = {Numerical Mathematics: Theory, Methods and Applications},
year = {2021},
volume = {14},
number = {4},
pages = {862--892},
issn = {2079-7338},
doi = {https://doi.org/10.4208/nmtma.OA-2021-0100},
url = {http://global-sci.org/intro/article_detail/nmtma/19522.html},
}

@article{liu2021iterated,
  title={Iterated numerical homogenization for multiscale elliptic equations with monotone nonlinearity},
  author={Liu, Xinliang and Chung, Eric and Zhang, Lei},
  journal={Multiscale Modeling \& Simulation},
  volume={19},
  number={4},
  pages={1601--1632},
  year={2021},
  publisher={SIAM}
}

@phdthesis{leitenmaier_Phdthesis,
  title={Analysis and numerical methods for multiscale problems in magnetization dynamics},
  author={Leitenmaier, Lena},
  year={2021},
  school={KTH Royal Institute of Technology}
}

@article{zhang_RPS,
  title={Polyharmonic homogenization, rough polyharmonic splines and sparse super-localization},
  author={Owhadi, Houman and Zhang, Lei and Berlyand, Leonid},
  journal={ESAIM: Mathematical Modelling and Numerical Analysis},
  volume={48},
  number={2},
  pages={517--552},
  year={2014},
  publisher={EDP Sciences}
}

@article{anrongbackward,
  title="{Analysis of backward Euler projection FEM for the Landau--Lifshitz equation}",
  author={An, Rong and Sun, Weiwei},
  journal={IMA Journal of Numerical Analysis},
  volume={42},
  number={3},
  pages={2336--2360},
  year={2022},
  publisher={Oxford University Press}
}

@article{AnewschemeLLG,
title = "{A new finite element scheme for Landau-Lifchitz equations}",
journal = {Discrete \& Continuous Dynamical Systems - S},
volume = {1},
number = {2},
pages = {187-196},
year = {2008},
author = {Alouges, Francois},
}

@article{santugini2007homogenization,
  title={Homogenization of ferromagnetic multilayers in the presence of surface energies},
  author={Santugini-Repiquet, K{\'e}vin},
  journal={ESAIM: Control, Optimisation and Calculus of Variations},
  volume={13},
  number={2},
  pages={305--330},
  year={2007},
  publisher={EDP Sciences}
}

@article{alouges2015homogenization,
  title={Homogenization of composite ferromagnetic materials},
  author={Alouges, Francois and Di Fratta, Giovanni},
  journal={Proceedings of the Royal Society A: Mathematical, Physical and Engineering Sciences},
  volume={471},
  number={2182},
  pages={20150365},
  year={2015},
  publisher={The Royal Society Publishing}
}

@article{choquet2018homogenization,
  title="{Homogenization of the Landau-Lifshitz-Gilbert equation in a contrasted composite medium}",
  author={Choquet, Catherine and Moumni, Mohammed and Tilioua, Mouhcine},
  journal={Discrete \& Continuous Dynamical Systems-S},
  volume={11},
  number={1},
  pages={35},
  year={2018},
  publisher={American Institute of Mathematical Sciences}
}

@article{alouges2021stochastic,
  title="{Stochastic homogenization of the Landau--Lifshitz--Gilbert equation}",
  author={Alouges, Francois and De Bouard, Anne and Merlet, Benoit and Nicolas, Lea},
  journal={Stochastics and Partial Differential Equations: Analysis and Computations},
  volume={9},
  pages={789--818},
  year={2021},
  publisher={Springer}
}

@article{UpscalingHMM,
author = {Leitenmaier, Lena and Runborg, Olof},
title = "{Upscaling Errors in Heterogeneous Multiscale Methods for the Landau--Lifshitz Equation}",
journal = {Multiscale Modeling \& Simulation},
volume = {20},
number = {1},
pages = {1-35},
year = {2022}
}

@article{leitenmaier2022heterogeneous,
  title="{Heterogeneous Multiscale Methods for the Landau--Lifshitz Equation}",
  author={Leitenmaier, Lena and Runborg, Olof},
  journal={Journal of Scientific Computing},
  volume={93},
  number={3},
  pages={76},
  year={2022},
  publisher={Springer}
}

@article{leitenmaier2022homogenization,
  title="{On homogenization of the Landau--Lifshitz equation with rapidly oscillating material coefficient}",
  author={Leitenmaier, Lena and Runborg, Olof},
  journal={Communications in Mathematical Sciences},
  volume={20},
  number={3},
  pages={653--694},
  year={2022},
  publisher={International Press of Boston}
}

@article{leitenmaier2023finite,
  title="{A finite element based Heterogeneous Multiscale Method for the Landau-Lifshitz equation}",
  author={Leitenmaier, Lena and Nazarov, Murtazo},
  journal={Journal of Computational Physics},
  volume={486},
  pages={112112},
  year={2023},
  publisher={Elsevier}
}

@article{chen2022multiscale,
  title="{On the Multiscale Landau--Lifshitz--Gilbert Equation: Two-Scale Convergence and Stability Analysis}",
  author={Chen, Jingrun and Du, Rui and Ma, Zetao and Sun, Zhiwei and Zhang, Lei},
  journal={Multiscale Modeling \& Simulation},
  volume={20},
  number={2},
  pages={835--856},
  year={2022},
  publisher={SIAM}
}

@article{owhadi2017multigrid,
  title={Multigrid with rough coefficients and multiresolution operator decomposition from hierarchical information games},
  author={Owhadi, Houman},
  journal={Siam Review},
  volume={59},
  number={1},
  pages={99--149},
  year={2017},
  publisher={SIAM}
}

@article{gutfleisch2011magnetic,
  title={Magnetic materials and devices for the 21st century: stronger, lighter, and more energy efficient},
  author={Gutfleisch, Oliver and Willard, Matthew A and Br{\"u}ck, Ekkes and Chen, Christina H and Sankar, SG and Liu, J Ping},
  journal={Advanced Materials},
  volume={23},
  number={7},
  pages={821--842},
  year={2011},
  publisher={Wiley Online Library}
}

@article{altmann2021numerical,
  title={Numerical homogenization beyond scale separation},
  author={Altmann, Robert and Henning, Patrick and Peterseim, Daniel},
  journal={Acta Numerica},
  volume={30},
  pages={1--86},
  year={2021},
  publisher={Cambridge University Press}
}

@article{henning2022superconvergence,
  title="{Superconvergence of time invariants for the Gross--Pitaevskii equation}",
  author={Henning, Patrick and W{\"a}rnegard, Johan},
  journal={Mathematics of Computation},
  volume={91},
  number={334},
  pages={509--555},
  year={2022}
}

@article{maier2022multiscale,
  title="{Multiscale scattering in nonlinear Kerr-type media}",
  author={Maier, Roland and Verf{\"u}rth, Barbara},
  journal={Mathematics of Computation},
  volume={91},
  number={336},
  pages={1655--1685},
  year={2022}
}

@article{verfurth2022numerical,
  title={Numerical homogenization for nonlinear strongly monotone problems},
  author={Verf{\"u}rth, Barbara},
  journal={IMA Journal of Numerical Analysis},
  volume={42},
  number={2},
  pages={1313--1338},
  year={2022},
  publisher={Oxford University Press}
}

@book{papanicolau1978asymptotic,
  title={Asymptotic analysis for periodic structures},
  author={Papanicolau, George and Bensoussan, Alain and Lions, J-L},
  year={1978},
  publisher={Elsevier}
}

@book{jikov2012homogenization,
  title={Homogenization of differential operators and integral functionals},
  author={Jikov, Vasili Vasilievitch and Kozlov, Sergei M and Oleinik, Olga Arsenievna},
  year={2012},
  publisher={Springer Science \& Business Media}
}

@article{dur91,
  title={Numerical calculation of equivalent grid block permeability tensors for heterogeneous porous media},
  author={Durlofsky, Louis J},
  journal={Water Resources Research},
  volume={27},
  number={5},
  pages={699--708},
  year={1991},
  publisher={Wiley Online Library}
}

@article{ab05,
  title={A multiscale finite element method for numerical homogenization},
  author={Allaire, Gr{\'e}goire and Brizzi, Robert},
  journal={Multiscale Modeling \& Simulation},
  volume={4},
  number={3},
  pages={790--812},
  year={2005},
  publisher={SIAM}
}

@article{weh02,
  title={Analysis of upscaling absolute permeability},
  author={Wu, Xiao-Hui and Efendiev, Yalchin and Hou, Thomas Y},
  journal={Discrete and Continuous Dynamical Systems Series B},
  volume={2},
  number={2},
  pages={185--204},
  year={2002},
  publisher={AIMS PRESS}
}

@article{ee03,
  title={The heterogeneous multiscale methods},
  author={E, Weinan and Engquist, Bjorn},
  journal={Communications in Mathematical Sciences},
  volume={1},
  number={1},
  pages={87--132},
  year={2003},
  publisher={International Press of Boston}
}

@article{abdulle2014analysis,
	Author = {Abdulle, Assyr and Vilmart, Gilles},
	Journal = {Mathematics of Computation},
	Number = {286},
	Pages = {513--536},
	Title = {Analysis of the finite element heterogeneous multiscale method for quasilinear elliptic homogenization problems},
	Volume = {83},
	Year = {2014}}

@article{ming2005analysis,
	Author = {E, Weinan and Ming, Pingbing and Zhang, Pingwen},
	Journal = {Journal of the American Mathematical Society},
	Number = {1},
	Pages = {121--156},
	Title = {Analysis of the heterogeneous multiscale method for elliptic homogenization problems},
	Volume = {18},
	Year = {2005}}

@book{berlyand2013introduction,
  title={Introduction to the network approximation method for materials modeling},
  author={Berlyand, Leonid and Kolpakov, Alexander G and Novikov, Alexei},
  number={148},
  series={Encyclopedia of Mathematics and its Applications},
  year={2013},
  publisher={Cambridge University Press}
}

@article{Arbogast_two_scale_04,
  title={Analysis of a two-scale, locally conservative subgrid upscaling for elliptic problems},
  author={Arbogast, Todd},
  journal={SIAM Journal on Numerical Analysis},
  volume={42},
  number={2},
  pages={576--598},
  year={2004},
  publisher={SIAM}
}

@article{egw10,
  title={Multiscale finite element methods for high-contrast problems using local spectral basis functions},
  author={Efendiev, Yalchin and Galvis, Juan and Wu, Xiao-Hui},
  journal={Journal of Computational Physics},
  volume={230},
  number={4},
  pages={937--955},
  year={2011},
  publisher={Elsevier}
}

@book{eh09,
  title={Multiscale finite element methods: theory and applications},
  author={Efendiev, Yalchin and Hou, Thomas Y},
  volume={4},
  year={2009},
  publisher={Springer Science \& Business Media}
}

@article{hughes98,
  title={The variational multiscale method—a paradigm for computational mechanics},
  author={Hughes, Thomas JR and Feij{\'o}o, Gonzalo R and Mazzei, Luca and Quincy, Jean-Baptiste},
  journal={Computer Methods in Applied Mechanics and Engineering},
  volume={166},
  number={1-2},
  pages={3--24},
  year={1998},
  publisher={Elsevier}
}

@article{bazilevs2007variational,
	Author = {Bazilevs, Y and Calo, VM and Cottrell, JA and Hughes, TJR and Reali, A and Scovazzi, G},
	Journal = {Computer Methods in Applied Mechanics and Engineering},
	Number = {1},
	Pages = {173--201},
	Publisher = {Elsevier},
	Title = {Variational multiscale residual-based turbulence modeling for large eddy simulation of incompressible flows},
	Volume = {197},
	Year = {2007}}

@article{berlyand2010flux,
	Author = {Berlyand, Leonid and Owhadi, Houman},
	Journal = {Archive for Rational Mechanics and Analysis},
	Number = {2},
	Pages = {677--721},
	Publisher = {Springer},
	Title = {Flux norm approach to finite dimensional homogenization approximations with non-separated scales and high contrast},
	Volume = {198},
	Year = {2010}}

@article{Owhadi:2011,
  title={Localized bases for finite-dimensional homogenization approximations with nonseparated scales and high contrast},
  author={Owhadi, Houman and Zhang, Lei},
  journal={Multiscale Modeling \& Simulation},
  volume={9},
  number={4},
  pages={1373--1398},
  year={2011},
  publisher={SIAM}
}

@article{owhadi2017gamblets,
  title={Gamblets for opening the complexity-bottleneck of implicit schemes for hyperbolic and parabolic ODEs/PDEs with rough coefficients},
  author={Owhadi, Houman and Zhang, Lei},
  journal={Journal of Computational Physics},
  volume={347},
  pages={99--128},
  year={2017},
  publisher={Elsevier}
}

@article{MalPet:2014,
  title={Localization of elliptic multiscale problems},
  author={M\"{a}lqvist, Axel and Peterseim, Daniel},
  journal={Mathematics of Computation},
  volume={83},
  number={290},
  pages={2583--2603},
  year={2014}
}

@article{egh12,
  title="{Generalized multiscale finite element methods (GMsFEM)}",
  author={Efendiev, Yalchin and Galvis, Juan and Hou, Thomas Y},
  journal={Journal of Computational Physics},
  volume={251},
  pages={116--135},
  year={2013},
  publisher={Elsevier}
}

@article{chung2014adaptiveDG,
  title={An adaptive generalized multiscale discontinuous Galerkin method for high-contrast flow problems},
  author={Chung, Eric T and Efendiev, Yalchin and Leung, Wing Tat},
  journal={Multiscale Modeling \& Simulation},
  volume={16},
  number={3},
  pages={1227--1257},
  year={2018},
  publisher={SIAM}}

@article{chung2015residual,
	title={Residual-driven online generalized multiscale finite element methods},
  author={Chung, Eric T and Efendiev, Yalchin and Leung, Wing Tat},
  journal={Journal of Computational Physics},
  volume={302},
  pages={176--190},
  year={2015},
  publisher={Elsevier}}

@article{chung2014adaptive,
  title="{An adaptive GMsFEM for high-contrast flow problems}",
  author={Chung, Eric T and Efendiev, Yalchin and Li, Guanglian},
  journal={Journal of Computational Physics},
  volume={273},
  pages={54--76},
  year={2014},
  publisher={Elsevier}
}

@article{ma2022novel,
  title={Novel design and analysis of generalized finite element methods based on locally optimal spectral approximations},
  author={Ma, Chupeng and Scheichl, Robert and Dodwell, Tim},
  journal={SIAM Journal on Numerical Analysis},
  volume={60},
  number={1},
  pages={244--273},
  year={2022},
  publisher={SIAM}
}

@article{di2020linear,
  title="{Linear second-order IMEX-type integrator for the (eddy current) Landau--Lifshitz--Gilbert equation}",
  author={Di Fratta, Giovanni and Pfeiler, Carl-Martin and Praetorius, Dirk and Ruggeri, Michele and Stiftner, Bernhard},
  journal={IMA Journal of Numerical Analysis},
  volume={40},
  number={4},
  pages={2802--2838},
  year={2020},
  publisher={Oxford University Press}
}

@article{praetorius2018convergence,
  title={Convergence of an implicit--explicit midpoint scheme for computational micromagnetics},
  author={Praetorius, Dirk and Ruggeri, Michele and Stiftner, Bernhard},
  journal={Computers \& Mathematics with Applications},
  volume={75},
  number={5},
  pages={1719--1738},
  year={2018},
  publisher={Elsevier}
}

@article{xie2020second,
  title={Second-order semi-implicit projection methods for micromagnetics simulations},
  author={Xie, Changjian and Garc{\'\i}a-Cervera, Carlos J and Wang, Cheng and Zhou, Zhennan and Chen, Jingrun},
  journal={Journal of Computational Physics},
  volume={404},
  pages={109104},
  year={2020},
  publisher={Elsevier}
}

@article{cai2022second,
  title="{A second-order numerical method for Landau-Lifshitz-Gilbert equation with large damping parameters}",
  author={Cai, Yongyong and Chen, Jingrun and Wang, Cheng and Xie, Changjian},
  journal={Journal of Computational Physics},
  volume={451},
  pages={110831},
  year={2022},
  publisher={Elsevier}
}

@article{gallistl2018numerical,
  title={Numerical Homogenization of H(curl)-Problems},
  author={Gallistl, Dietmar and Henning, Patrick and Verf{\"u}rth, Barbara},
  journal={SIAM Journal on Numerical Analysis},
  volume={56},
  number={3},
  pages={1570--1596},
  year={2018},
  publisher={SIAM}
}

@article{henning2020computational,
  title="{Computational homogenization of time-harmonic Maxwell's equations}",
  author={Henning, Patrick and Persson, Anna},
  journal={SIAM Journal on Scientific Computing},
  volume={42},
  number={3},
  pages={B581--B607},
  year={2020},
  publisher={SIAM}
}

@article{freese2024computational,
  title={Computational multiscale methods for nondivergence-form elliptic partial differential equations},
  author={Freese, Philip and Gallistl, Dietmar and Peterseim, Daniel and Sprekeler, Timo},
  journal={Computational Methods in Applied Mathematics},
  volume={24},
  number={3},
  pages={649--672},
  year={2024},
  publisher={De Gruyter}
}

@book{chen1998second,
  title={Second order elliptic equations and elliptic systems},
  author={Chen, Ya-Zhe and Wu, Lan-Cheng},
  volume={174},
  year={1998},
  publisher={American Mathematical Soc.}
}

@article{doding2022two,
  title="{A two level approach for simulating Bose-Einstein condensates by localized orthogonal decomposition}",
  author={D{\"o}ding, Christian and Henning, Patrick and W{\"a}rnegard, Johan},
  journal={arXiv preprint arXiv:2212.07392},
  year={2022}
}

@article{bonizzoni2022super,
  title={Super-localized orthogonal decomposition for convection-dominated diffusion problems},
  author={Bonizzoni, Francesca and Freese, Philip and Peterseim, Daniel},
  journal={BIT Numerical Mathematics},
  volume={64},
  number={3},
  pages={33},
  year={2024},
  publisher={Springer}
}

@article{bonizzoni2024reduced,
  title={A reduced basis super-localized orthogonal decomposition for reaction-convection-diffusion problems},
  author={Bonizzoni, Francesca and Hauck, Moritz and Peterseim, Daniel},
  journal={Journal of Computational Physics},
  volume={499},
  pages={112698},
  year={2024},
  publisher={Elsevier}
}

\section{Appendix}
\label{sec:appendix}

\subsection{Various Time Discretization Schemes}
\label{subsec:Various Time Discretization Schemes}

\paragraph{Prohl's method} The fully implicit scheme proposed by Prohl \cite{prohl2001computational} ensures numerical stability through complete nonlinear treatment:
\begin{equation}
    \frac{1}{\Delta t}(\m_h^{n+1} - \m_h^{n},\v_h) - \lambda( \h_{\eff}^{n+1},\v_h) + ( \m_h^{n+1} \times \h_{\eff}^{n+1} , \v_h) = -\lambda ((\m_h^{n} \cdot \h_{\eff}^{n}) \cdot \m_h^{n+1} , \v_h),
    \label{eqn:fullynonlinear}
\end{equation}
for all $\v_h \in \V_h$ and $n=0,1,\cdots,N-1$. This formulation requires nonlinear solvers at each time step.

\paragraph{Cimr\'ak's method} 
Cimr\'ak \cite{cimrak2005error} developed a semi-implicit variant by linearizing the cross product term:
\begin{equation}
    \frac{1}{\Delta t}(\m_h^{n+1} - \m_h^{n},\v_h) - \lambda( \h_{\eff}^{n+1},\v_h) + ( \m_h^{n} \times \h_{\eff}^{n+1} , \v_h) = -\lambda ( (\m^{n} \cdot \h_{\eff}^{n}) \cdot \m_h^{n+1} , \v_h).
    \label{Cimark}
\end{equation}
This modification yields a linear system solvable without iteration, significantly improving computational efficiency while maintaining stability.

\paragraph{Gao's method} 
Gao \cite{gao2014optimal} further simplified the right-hand side treatment:
\begin{equation}
    \frac{1}{\Delta t}(\m_h^{n+1} - \m_h^{n},\v_h) - \lambda( \h_{\eff}^{n+1},\v_h) + ( \m_h^{n} \times \h_{\eff}^{n+1} , \v_h) = -\lambda ( (\m^{n} \cdot \h_{\eff}^{n}) \cdot \m_h^{n} , \v_h).
    \label{GaoFEM}
\end{equation}
This scheme achieves optimal $L^2$ and $H^1$ error estimates without restrictive time-step conditions ($\Delta t=O(h^\alpha)$). The unit length constraint relaxation yields:
\begin{equation}
    ||1 - |\m_h^n|^2 ||_{L^2} \le C_0 (\Delta t + h^{2}),
\end{equation}
demonstrating controlled deviation from the unit sphere.

\paragraph{An's method} 
An \cite{anrongbackward} introduced a projection-based approach:
\begin{equation}
    \frac{1}{\Delta t}(\tilde{\m}_h^{n+1} - \m_h^{n},\v_h) - \lambda( \h_{\eff}^{n+1},\v_h) + ( \m_h^{n} \times \h_{\eff}^{n+1} , \v_h) = -\lambda ( (\m_h^{n} \cdot \Bar{\h}_{\eff}^{n+1}) \cdot \m_h^{n} , \v_h)
    \label{AnFEM}
\end{equation}
with post-processing projection step $\m_h^{n+1} = \tilde{\m}_h^{n+1}/|\tilde{\m}_h^{n+1}|$, and $\Bar{\h}_{\eff}^{n+1} = \nabla  \tilde{\m}_h^{n+1} \cdot \kappa \nabla \m_h^n + \lambda \tilde{\m}_h^{n+1} \cdot \m_{h,a}^n$. While theoretically requiring $\Delta t=O(\eps_0 h)$ and $\m_0 \in H^{r+1}$ ($r\geq2$), numerical experiments show these conditions can be relaxed.

\subsection{Measurement Function Types}
\label{subsec:Two types of measurement functions}
Let the set of edge-based measurement functions be denoted by $\Phi_{\mathcal{E}}=\{ \phi_{e} \}_{e \in \mathcal{E}_H}$ and the set of volume-based measurement functions by $\Phi_{\mathcal{T}}=\{ \phi_{\tau} \}_{\tau \in \mathcal{T}_{H}}$.
The GRPS framework mainly employs two primary measurement function types:

\begin{itemize}
    \item \textbf{Edge-based (Case E)}: For $e \in \mathcal{E}_H$,
    \begin{equation}
        \phi_{e} = |e|^{\frac{2-d}{2(d-1)}} \chi(e), \quad \langle \phi_{e}, u \rangle = \int_{e} u \d S;
    \end{equation}
    
    \item \textbf{Volume-based (Case V)}: For $\tau \in \mathcal{T}_{H}$,
    \begin{equation}
        \phi_{\tau} = \sqrt{|\tau|} \chi(\tau), \quad \langle \phi_{\tau}, u \rangle = \sqrt{|\tau|}\int_{\tau} u \d \x,
    \end{equation}
\end{itemize}
where $\chi$ denotes characteristic functions and $|\cdot|$ the geometric measure.

\subsection{The multiscale trigonometric (MsTrig) coefficient}
\label{subsec:six-scale multiscale coefficient}
The functions in \eqref{multiscale_coefficient_k} is defined by
\begin{itemize}
    \item $f_1(x,y,\eps_1) = \frac{1.1+\sin(2\pi x/\eps_1)}{1.1+\sin(2\pi y/\eps_1)}$,
    \item $f_2(x,y,\eps_2) = \frac{1.1+\sin(2\pi y/\eps_2)}{1.1+\cos(2 \pi x/\eps_2)}$,
    \item $f_3(x,y,\eps_3) = \frac{1.1+\cos(2\pi x/\eps_3)}{1.1+\sin(2\pi y/\eps_3)}$,
    \item $f_4(x,y,\eps_4) = \frac{1.1+\sin(2\pi y/\eps_4)}{1.1+\cos(2\pi x/\eps_4)}$,
    \item $f_5(x,y,\eps_5) = \frac{1.1+\cos(2\pi x/\eps_5)}{1.1+\sin(2\pi y/\eps_5)}$,
\end{itemize}
where $\eps_1 = 1/5$, $\eps_2 = 1/13$, $\eps_3 = 1/17$, $\eps_4 = 1/31$, $\eps_5 = 1/65$. 





\end{document}